\documentclass[11pt]{article} 

\usepackage[utf8]{inputenc} 

\usepackage{geometry} 
\usepackage{graphicx} 
\usepackage{rotating}

\usepackage{booktabs} 
\usepackage{array} 
\usepackage{paralist} 
\usepackage{verbatim} 
\usepackage{subfig} 
\usepackage{amsmath,amssymb,verbatim}
\usepackage{color}
\usepackage{graphicx}
\usepackage{mathtools}
\usepackage{algorithm}
\usepackage{algpseudocode}
\usepackage{bigfoot}
\usepackage{authblk}

\mathtoolsset{showonlyrefs}

\newcommand{\R}{\mathbb{R}}
\newcommand{\Rad}{\mathcal{R}}

\newcommand{\betah}{\frac{\beta}{2}}
\newcommand{\Radlp}{\mathcal{R}_{\mathrm{lp}}}
\newcommand{\Radlppart}{{\mathcal{R}_\mathrm{lp}^{\mathrm{p}}}}
\newcommand{\midx}{m}

\newcommand{\Jlp}{\mathcal{W}_{\mathrm{lp}}}
\newcommand{\thetalp}{{\theta_{\mathrm{lp}}}}
\newcommand{\slp}{{s_{\mathrm{lp}}}}
\newcommand{\thetap}{{\theta_\mathrm{p}}}
\DeclareMathOperator{\sinc}{sinc}

\graphicspath{{./figs/}}
\title{Fast algorithms and efficient GPU implementations for the Radon transform and the back-projection operator represented as convolution operators}

\author{Fredrik Andersson} 
\author{Marcus Carlsson} 
\author{Viktor V Nikitin} 
\affil{Center for Mathematical Sciences, Lund University, Box 118, 22100 Lund, Sweden}
\date{} 

\begin{document}
\maketitle


\begin{abstract}
The Radon transform and its adjoint, the back-projection operator, can both be expressed as convolutions in log-polar coordinates. Hence, fast algorithms for the application of the operators can be constructed by using FFT, if data is resampled at log-polar coordinates. Radon data is typically measured on an equally spaced grid in polar coordinates, and reconstructions are represented (as images) in Cartesian coordinates. Therefore, in addition to FFT, several steps of interpolation have to be conducted in order to apply the Radon transform and the back-projection operator by means of convolutions. However, in comparison to the interpolation conducted in Fourier-based gridding methods, the interpolation performed in the Radon and image domains will typically deal with functions that are substantially less oscillatory. Reasonable reconstruction results can thus be expected using interpolation schemes of moderate order. It also provides better control over the artifacts that can appear due to measurement errors. 

Both the interpolation and the FFT operations can be efficiently implemented on Graphical Processor Units (GPUs). For the interpolation, it is possible to make use of the fact that linear interpolation is hard-wired on GPUs, meaning that it has the same computational cost as direct memory access. Cubic order interpolation schemes can be constructed by combining linear interpolation steps which provides important computation speedup. 

We provide details about how the Radon transform and the back-projection can be implemented efficiently as convolution operators on GPUs. For large data sizes, speedups of about 10 times are obtained in relation to the computational times of other software packages based on GPU implementations of the Radon transform and the back-projection operator. Moreover, speedups of more than a 1000 times are obtained against the CPU-implementations provided in the MATLAB image processing toolbox.
\end{abstract}

\section{Introduction}
The two-dimensional Radon transform is the mapping of functions to their line integrals, i.e., a mapping $\Rad: \R^2\to S^1\times \R$ where $S^1$ denotes the unit circle, defined by
\begin{align}
\Rad f(\theta,s)=\int f(x)\delta(x\cdot\theta-s)d x.
\label{radon_def}
\end{align}
The parameter $\theta$ represents the (normal) direction of the lines, and the parameter $s$ denotes the (signed) distance of the line to the origin. It is customary to use $\theta$ both as a point on the unit sphere and as an angle, i.e., the notation $x\cdot\theta$ is used to parameterize lines in $\R^2$ by the relation $x\cdot\theta=x_1\cos(\theta)+x_2\sin(\theta)$. Note that each line is defined twice in this definition, since $s$ can take both positive and negative, and since $\theta \in S^1$.

\begin{figure}
\centering 	
\subfloat{\includegraphics[clip=true,trim=45mm 50mm 40mm 0mm,width=0.48\textwidth]{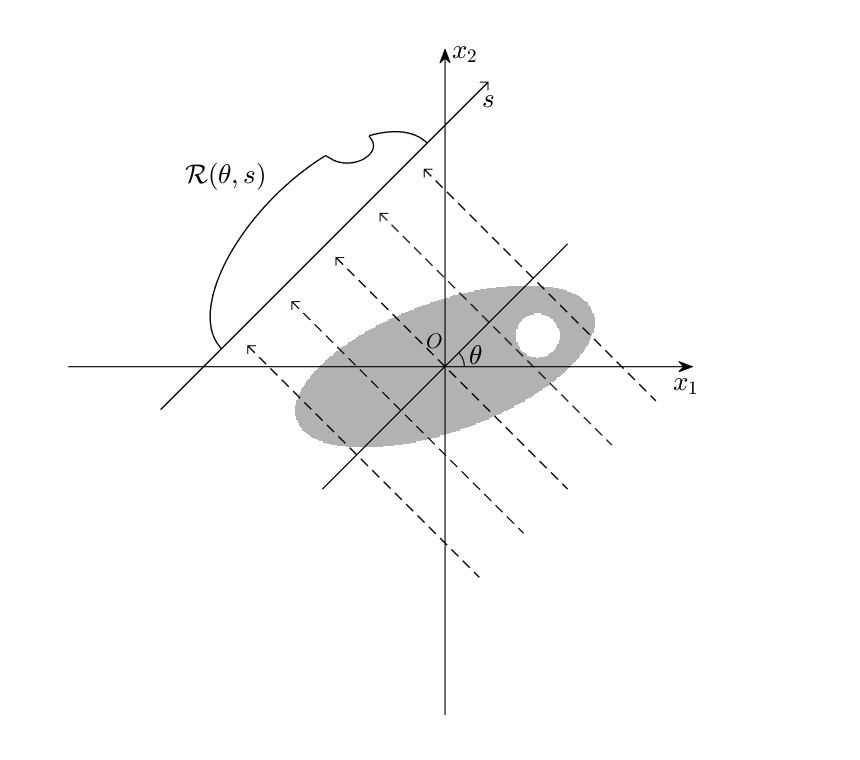}}
\caption{Scheme of computing projections for a given angle $\theta$.}
\label{fig:Rexpl}
\end{figure}

A schematic illustration of the Radon transform is depicted in Figure \ref{fig:Rexpl}, where beams are propagating through an object and after absorption are measured by receivers. The Radon transform appears for instance in computational tomography (CT). The tomographic inversion problem lies in recovering an unknown function $f$ given knowledge of $\Rad f$. For more details about CT, see \cite{faridani2003introduction,frikel2013characterization,kak1988principles,natterer1986computerized}.

One of the most popular methods of inverting the usual Radon transform is by means of the filtered back-projection (FBP) method \cite{natterer1986computerized}. It uses the inversion formula
\begin{equation}\label{FBP}
f = \Rad^\# \mathcal{W} \Rad f.
\end{equation}
where $\mathcal{W}$ is a convolution operator acting only on the $s$-variable, and where $\Rad^\# :S^1\times \R \to R^2$ denotes the back-projection operator, which integrates over all lines through a point, i.e.,
\begin{equation}
\Rad^\# g(x) = \int_{S^1} g(\theta,x\cdot \theta).
\end{equation}
The back-projection operator is adjoint to the Radon transform. The convolution operator $\mathcal{W}$ can either be describes as a Hilbert transform followed by a derivation, both with respect to the variable $s$; or as a convolution operator with a transfer function being a suitable scaled version of $|\sigma|$, where $\sigma$ denotes the conjugate variable of $s$.

The filtered back-projection algorithm has a  time complexity of $N^3$, if we assume that reconstructions are made on an $N\times N$ lattice and that the numbers of samples in $s$ and $\theta$ are both  $\mathcal{O}(N)$. This is because for each point $x$, integration has to be made over all lines ($N$ directions) passing through that point. However, there are methods for fast $\mathcal{O}(N^2\log N)$  back-projection, see for instance \cite{basu2000n, danielsson1997backprojection, george2007fast}.

\begin{figure}[t!]
\begin{minipage}{.48\linewidth}
\subfloat[][]{\includegraphics[width=0.49\linewidth,height=0.49\linewidth]{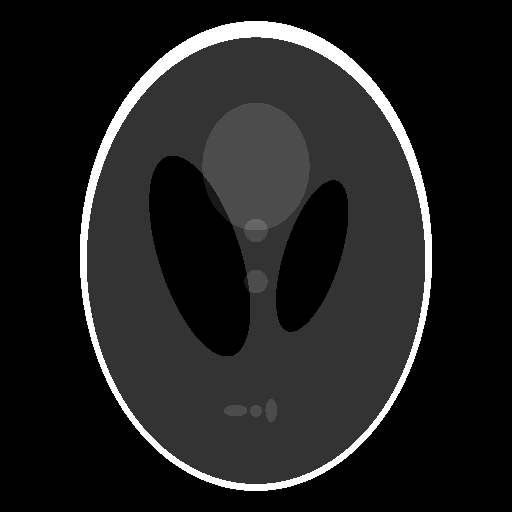}}
\subfloat[][]{\includegraphics[width=0.49\linewidth,height=0.49\linewidth]{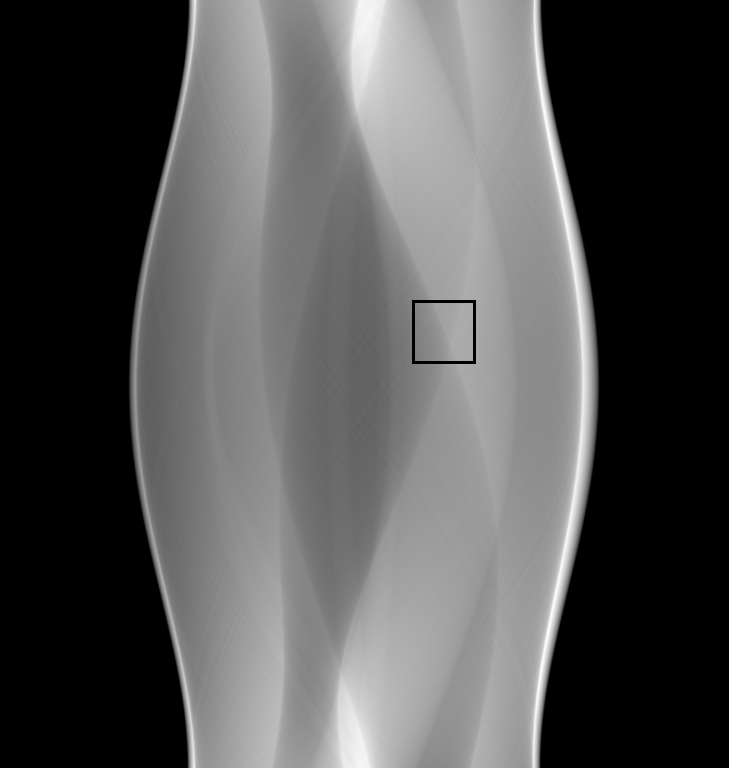}}
\end{minipage}
\subfloat[][]{
\begin{minipage}{.24\linewidth}
\centering
\hspace{-1.5cm} \includegraphics[width=0.4\linewidth,height=0.4\linewidth]{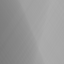} \\\vspace{0.5cm}\hspace{1.5cm}
\includegraphics[width=0.4\linewidth,height=0.4\linewidth]{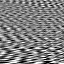}
\end{minipage}
}
\subfloat[][]{
\begin{minipage}{.24\linewidth}
\includegraphics[width=\linewidth,height=\linewidth]{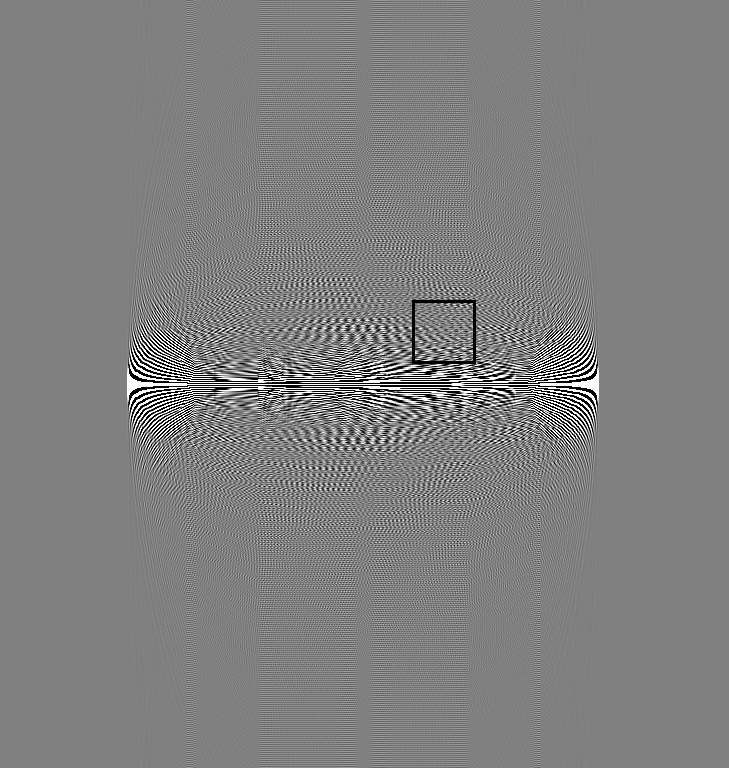}
\end{minipage}
}
\caption{\label{fig:phantom} The (modified) Shepp-Logan phantom (a) and its Radon transform (b). The panel (d) shows the one-dimensional Fourier transform of the Radon transform with respect to $s$. The regions indicated by the black squares in (b) and (d) are shown in higher resolution in (c).}
\end{figure}

Another class of fast methods inversion of Radon data goes via the Fourier-slice theorem. This is a result by which the Radon data can be mapped to a polar sampling of the unknown function in the frequency domain. To recover the unknown function, interpolation-like operations (for instances the ones presented in \cite{USFFT,brandt2000fast,kalamkar2012high}) have to be applied in the frequency domain. The data in the frequency domain will typically be oscillatory, as seen in Figure \ref{fig:phantom}d) and (at in increased resolution) in the lower right panel of Figure \ref{fig:phantom}c). Hence, in order to interpolate data of this type high interpolation order is required. In comparison, the data in the Radon domain will not be particularly oscillatory. This is illustrated in  Figure \ref{fig:phantom}b) and the upper left panel of Figure \ref{fig:phantom}c). This fact implies that interpolation methods of moderate order can be expected to produce reasonable results. This means in turn that less time can be spent on conducting interpolation. It also gives more control over the interpolation errors, as local errors will be kept local in the Radon domain (and hence more easily distinguishable as artifacts in the reconstruction).

In this paper, we discuss how to design fast algorithms for the application of the Radon transform and the back-projection operator by using the fact that they can be expressed in terms of convolutions when represented in log-polar coordinates \cite{andersson2005fast,eggermont1983tomographic,hansen1981theory,verly1981circular}. In particular, we follow the approach suggested in \cite{andersson2005fast}. This formulation turns out to be particularly well-suited for implementation on GPUs. A major advantage with using GPUs is that the routines for linear interpolation are fast. In fact, the cost for computing linear interpolation is the same as reading directly from memory \cite{govindaraju2006memory}. This feature can be utilized for constructing fast interpolation schemes. In particular, in this paper we will work with cubic interpolation on GPU \cite{ruijters2008efficient,sigg2005fast}.

For the sake of comparison, we will provide performance and accuracy tests of the proposed method along with comparisons against other software packages for tomographic computations. We also conduct a performance comparison between the different methods when utilized in iterative reconstruction techniques. The iterative methods rely on applying the Radon transform and the back-projection operator several times. An advantage of keeping all computations on the GPU is that the needed time for CPU-GPU memory transfer can be reduced.

\section{The Radon transform and the back-projection expressed as convolutions}\label{secrad}
We recapitulate some of the main ideas of the method described in \cite{andersson2005fast}. A key part there is the usage of log-polar coordinates, i.e.,
\begin{align}
\begin{cases}
&x_1=e^\rho\cos(\theta),\\
&x_2=e^\rho\sin(\theta),
\end{cases}
\label{param}
\end{align}
where $-\pi<\theta<\pi$. To simplify the presentation, we identify $f(\theta,\rho)$ with $f(x_1,x_2)$ if $(\theta,\rho)$ in the log-polar coordinate system corresponds to the point $x=(x_1,x_2)$ the Cartesian coordinate system, and similarly for other coordinate transformations.

By representing the distance between lines and the origin $s=e^\rho$, and by a change of variables in \eqref{radon_def} from Cartesian to log-polar coordinates 
the log-polar Radon transform  can be expressed  as
\begin{align}
\Radlp f(\theta,\rho)&=\int_{-\pi}^{\pi}\int_{-\infty}^{\infty}f(\theta',\rho')e^{\rho'}\delta\left(\cos(\theta-\theta')-e^{\rho-\rho'}\right) d\rho'd \theta' \nonumber \\
&=\int_{-\pi}^{\pi}\int_{-\infty}^{\infty}f(\theta',\rho')e^{\rho'}\zeta \left( \theta-\theta',\rho-\rho'\right) d\rho' d \theta' .
\label{logdef}
\end{align}
In particular, $\Rad f(\theta,s) = \Radlp f(\theta,\log(s))$ for $s>0$ (which is sufficient since the information for $s<0$ is redundant).

We briefly mention how to put this formula in a theoretical framework. Set $S=(-\pi,\pi)\times \mathbb{R}$ and note that, for a compactly supported smooth function $h$ on $S$, we have $$\int_{-\pi}^{\pi}\int_{-\infty}^{\infty}h(\theta,\rho)\zeta \left( \theta,\rho\right) d\rho d \theta = \int_{-\pi/2}^{\pi/2}\frac{h(\theta,\log\left( \cos(\theta)\right))}{\cos\theta} d \theta,$$ which can be written as $\int_S h d\mu$ where $\mu$ is an infinite measure on $S$. Hence, the formula extends by continuity to, e.g., all continuous compactly supported functions $h$ in $S$. It follows that \eqref{logdef} is well-defined whenever $f$ is a continuous function which is zero in a neighborhood of 0.

As the Radon transform in the coordinate system $(\theta,\rho)$ is essentially a convolution between $f$ and the distribution $\zeta(\theta,\rho) = \delta( \cos(\theta)-e^{\rho})$, it can be rapidly computed by means of Fourier transforms. Special care has to be taken to the distribution $\zeta$, an issue we will return to in what follows. Ignoring possible difficulties with the distribution $\zeta$, let us discuss how \eqref{logdef} can be realized by using fast Fourier transforms. It is natural to assume that the function $f$ has compact support (the object that is measured has to fit in the device that is measuring it). The compact support also implies that the Radon transform of $f$ will have compact support in the $s$-variable. However, this is not true in the log-polar setting, since $\rho \rightarrow -\infty$ as $s \rightarrow 0$. 

\begin{figure}
\centering 	
\subfloat[][]{\includegraphics[trim = 45mm 35mm 25mm 35mm,clip=true,width=0.48\textwidth]{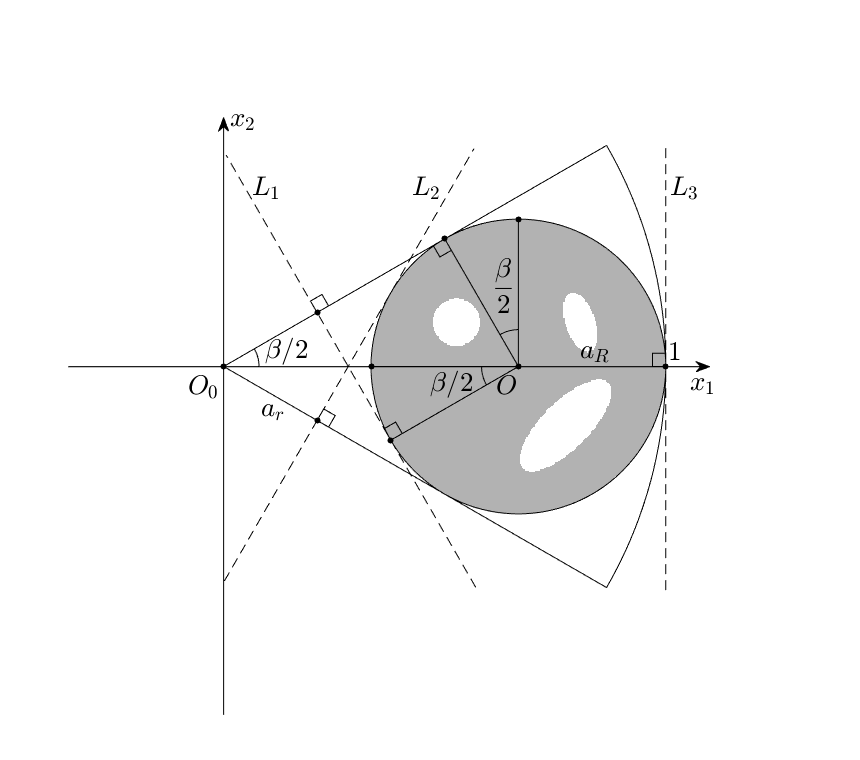}}
\subfloat[][]{\includegraphics[trim = 25mm 25mm 25mm 25mm,clip=true,width=0.48\textwidth]{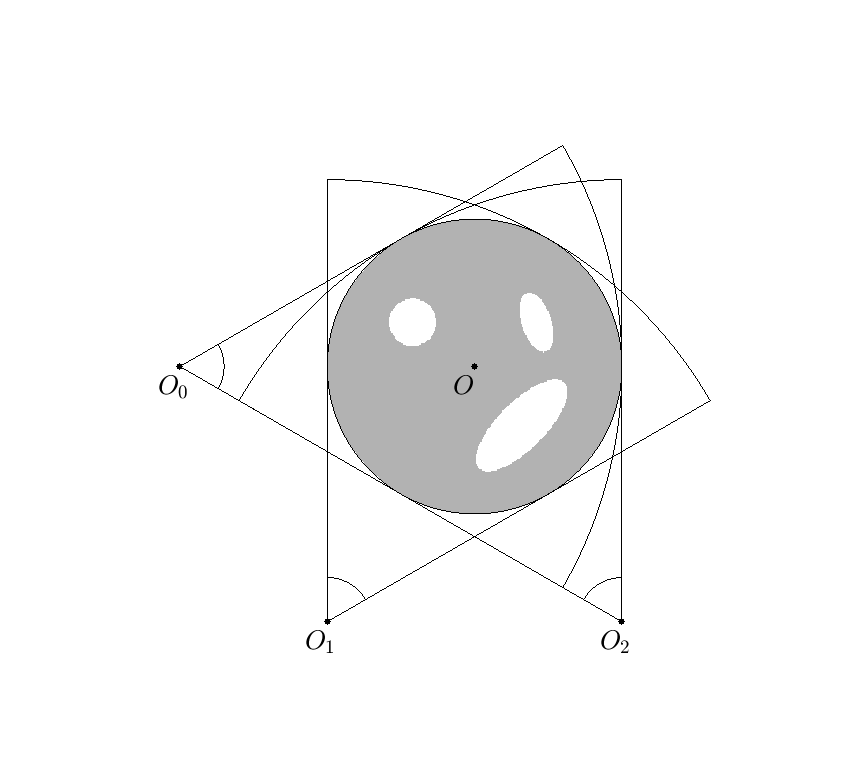}} \vspace{0.01mm} \caption{(a) Tangent lines to the circle to determine the support of the log-polar Radon transform function; (b) three angle spans to compute partial Radon transforms.}
\label{fig:spans}
\end{figure}

Note also that for any point $x$ in the plane there is a direction for which there is a line passing through $x$ and the origin. This implies that it is not possible to approximate the values of the Radon transform by using a finite convolution in log-polar coordinates if it is to be computed for all possible line directions. However, by restricting the values of $\theta$, and by making a translation so that the support of $f$ is moved away from the origin, it is in fact possible to describe the \emph{partial Radon transform} as a finite convolution, and then recover the full Radon transform by adding the contributions from various partial Radon transforms. The setup is illustrated in Fig. \ref{fig:spans}.

More precisely, let $\beta$ be a fixed angle and let $a_R$, $a_r,$ $L_1,~L_2$ and $L_3$ be as in Figure \ref{fig:spans}a). From the geometry we see that
\begin{align}
&a_R=\frac {\sin\left(\betah\right)}{1+\sin\left(\betah\right)}, \quad \mbox{ and } \quad a_r=\frac{\cos\left(\betah\right)-\sin\left(\betah\right)}{1+\sin\left(\betah\right)}.
\label{aconst}
\end{align}
Assume for the moment that $f$ has support in the gray circle indicated in Figure \ref{fig:spans}a). In log-polar coordinates $(\theta,\rho)$ it then has support inside of $\left[\log (1-2a_R),0\right]\times\left[-\betah,\betah\right]$. If we restrict our attention to values of $\Radlp(f)$ in the same sector $\theta\in \left[-\betah,\betah\right]$, then the only nonzero values of $\Radlp(f)$ will be for $\rho$ in the interval $\left[\log a_r,0\right]$.
This means that $\Radlp(f)(\theta,\rho)$ for these values can be computed by the finite convolution
\begin{align}
&\int_{-\beta/2}^{\beta/2}\int_{\log (1-2a_R) }^{0}  f(\theta',\rho') e^{\rho'}\zeta \left( \theta-\theta',\rho-\rho'\right)  d \theta' d\rho',
\label{logdefsimple}
\end{align}
where $(\theta,\rho)\in\left[-\betah,\betah\right]\times \left[\log a_r,0 \right].$ We now replace the integral \eqref{logdefsimple} by the periodic convolution
\begin{align}
&\Radlppart f(\theta,\rho)=\int_{-\beta}^{\beta}\int_{\log (a_r) }^{0}  f(\theta',\rho') e^{\rho'} \zeta_{\mathrm{per}} \left( \theta-\theta',\rho-\rho'\right)  d \theta' d\rho',
\label{lograd_periodic}
\end{align}
where $\zeta_{\mathrm{per}}$ is the periodic extension of $\zeta$ defined on $[\log(a_r),0]\times [-\beta,\beta]$. It is readily verified that for $(\theta,\rho)\in\left[-\betah,\betah\right]\times \left[\log a_r,0 \right]$, it thus holds that
\begin{align*}
\Radlppart f(\theta,\rho)&=\Rad f(\theta,e^\rho).
\end{align*}
We refer to  $\Radlppart$ as the \emph{partial log-polar Radon transform}.

Note that, in analogy with the argument following \eqref{logdef}, the formula \eqref{logdefsimple} can be written as a convolution between $f(\theta',\rho') e^{\rho'}$ and a finite measure, whereas \eqref{lograd_periodic} can be written as a convolution with a locally finite periodic measure. The above formulas are thus well defined as long as $f$ is continuous (or even piecewise continuous) in the log-polar coordinates. We refer to \cite[Chapter 11]{zemanian}, for basic results about convolution between functions and periodic distributions.

\begin{figure}
\centering 	
\subfloat[][]{\includegraphics[trim = 2cm 2cm 2cm 2cm,width=0.48\textwidth]{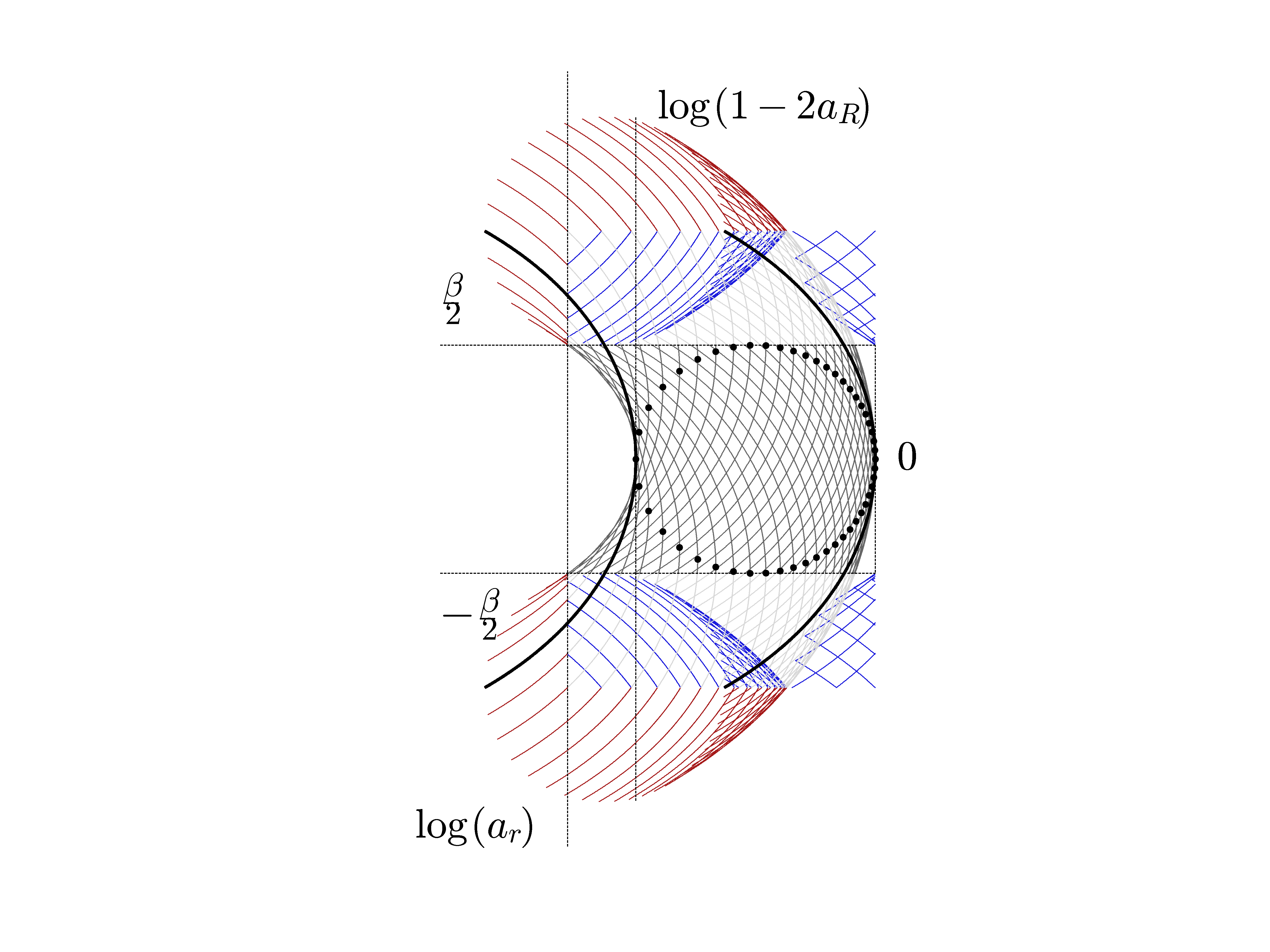}}
\subfloat[][]{\includegraphics[trim = 2cm 2cm 2cm 2cm,width=0.48\textwidth]{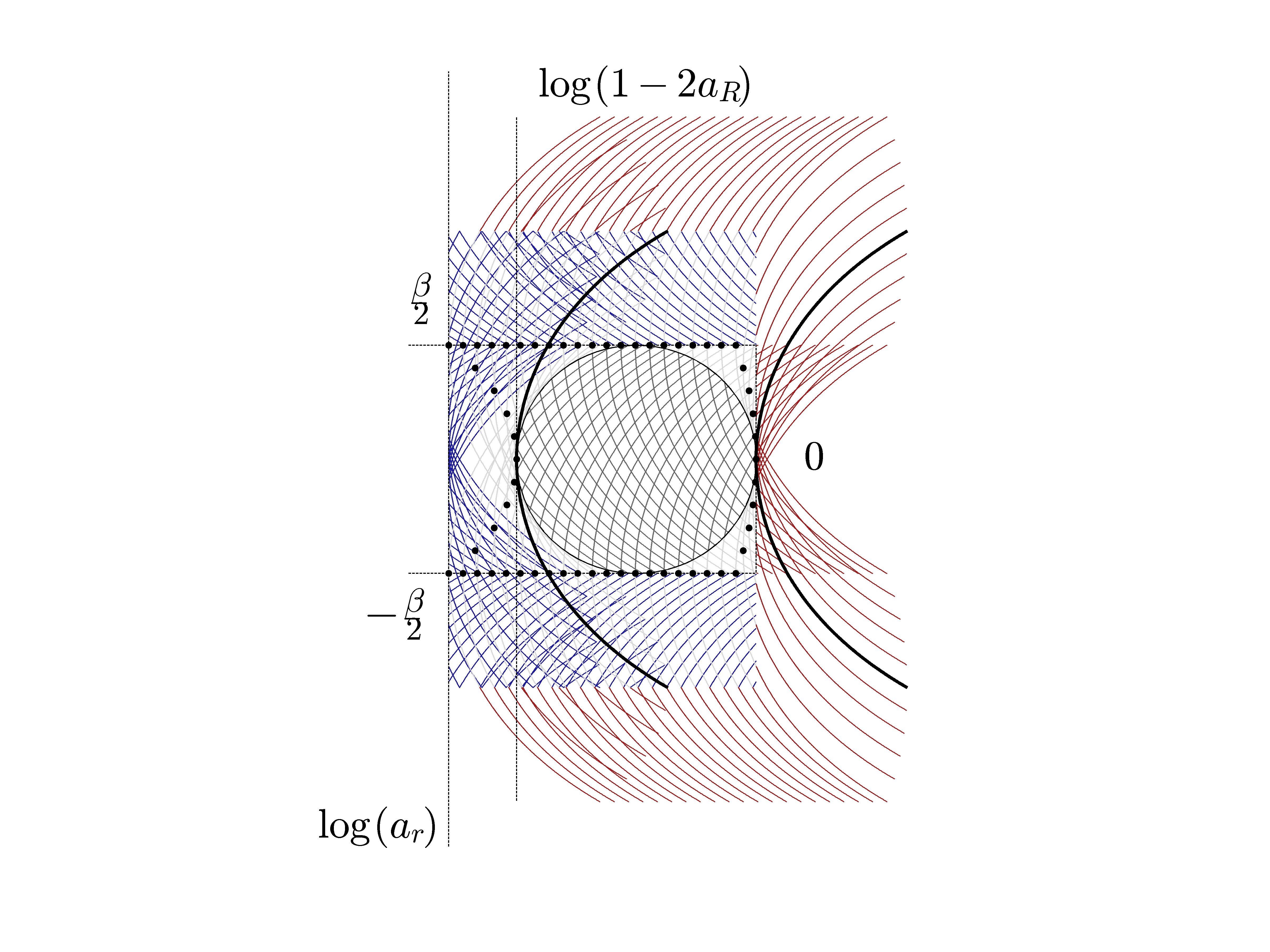}}
\caption{\label{fig:conv_support} The effect on the support due to convolutions: a) Radon transform; b) Back-projection.}
\end{figure}
The convolution setup for the Radon transform is depicted in Figure \ref{fig:conv_support}a). The rightmost black solid curve $C$ shows $\rho=\log(\cos(\theta))$, ($-\beta<\theta<\beta$), which is the support of $\zeta$ in $\left[-\beta,\beta\right]\times \left[\log a_r,0 \right]$. Let $D$ denote the circle
\begin{equation}\label{circle_D}
D=\{ (x_1,x_2): (x_1-1+a_R)^2+x_2^2<a_R^2\}
\end{equation}
in log-polar coordinates. The black dots shows the perimeter of $D$, and the gray curves indicate translation of the curve $C$ associated with the black dots on the circle, within the interval $[-\beta,\beta]$. There is a difference in grayscale for the points with $\theta$ inside the range $[-\betah,\betah]$, as only these values are of interest to us. Note that the smallest $\rho$-value of the contributing part in this interval is $\log(a_r)$. Moreover, the red lines indicate parts of the translations of $C$ outside $[\log(a_r),0]\times[-\beta,\beta]$, and the blue lines shows how these curves are wrapped back into the domain $[\log(a_r),0]\times [-\beta,\beta]$ by the periodic extension of $\zeta$. We see that these alias effects do not have any influence on the domain $[\log(a_r),0]\times [-\betah,\betah]$.

We now describe how $\Radlp^{\beta}$ can be used to recover $\Rad f$ for a function $f$ with support in the unit circle.
We split the angular variable into $M$ different parts, $\beta=\frac{\pi}{M}$. For $m=0,1,\ldots M-1$ let $\mathsf{T}_m : \mathbb{R}^2\rightarrow \mathbb{R}^2$ denote the change of coordinates
\begin{equation} \label{Tdef}
\mathsf{T}_m(x) =
a_R\begin{pmatrix}
\cos(\midx\beta) & \sin(\midx\beta)\\
-\sin(\midx\beta) & \cos(\midx\beta)
\end{pmatrix}
\begin{pmatrix}
x_1\\
x_2
\end{pmatrix}+
\begin{pmatrix}
1-a_R\\
0
\end{pmatrix} ,
\end{equation}
and let $T_m f = f(\mathsf{T}_m^{-1} x)$. Note that \begin{equation}\label{Radon_tm}\Rad f (\theta,s)=a_R^{-1}\Rad \left(T_{m} f\right) \Big(\theta-m\beta,a_Rs+(1-a_R)\cos(\theta-m\beta)\Big).\end{equation}
The Radon transform can thus be computed for arbitrary $\theta$ and $0<s<1$, by using the relation
\begin{equation}\label{Radon_all_pb}
\Rad f(\theta,s) = a_R^{-1} \Radlppart \left(T_{m} f\right)\Big(\theta-m\beta,\log\left(a_R s +(1-a_R)\cos\left(\theta-m\beta\right)\right)\Big),
\end{equation}
where $m= [\theta/\beta (\bmod M)]$ and $[x]$ denotes the rounding operator to the closest integer to $x$, and where $\bmod$ denotes the modulus operator.

We denote the change of coordinates above by
\begin{equation}\label{Sm_def}
\mathsf{S}_m(\theta,s) = \Big( \theta-m\beta,\log\left( a_Rs+(1-a_R)\cos(\theta-m\beta) \right) \Big),
\end{equation}
and we then have that
\begin{equation}
\Rad f(\theta,s)=a_R^{-1} \Radlppart \left(T_{m} f\right)(\mathsf{S}_m(\theta,s)) .
\end{equation}

We remark that for fixed $\theta$, the connection \eqref{Radon_tm} between the Radon data in the different domains has an affine dependence on $s$. Since the filter operator $\mathcal{W}$ acts as a convolution operator with regards to $s$, its action will in principle be the same regardless if the coordinate transformation $T_m$ is used or not. We use the notation $\Jlp$ to denote the action of the filter operator in log-polar coordinates.

%
The adjoint operator (back-projection) associated with the Radon transform \eqref{radon_def} can be written as
\begin{align}
\Rad^\# g(x) = \int_{-\infty}^{\infty}\int_{S^1} g(\theta,s)\delta(x\cdot \theta - s) \, d\theta ds,
\label{def_Radast}
\end{align}
cf. \cite{natterer1986computerized}. It is a weighted integral of $g$ over lines passing through the point $x$, and just as for the Radon transform it can be expressed as a convolution in log-polar coordinates. We define
\begin{align*}
\Radlp^\# g(\theta,\rho) &= \int_{-\infty}^{\infty} \int_{-\pi}^\pi g(\theta',\rho') \delta\left(e^{\rho-\rho'} \cos(\theta-\theta')-1\right) \, d\theta' d\rho' \\
&= \int_{-\infty}^{\infty} \int_{-\pi}^\pi g(\theta',\rho') \zeta^\# \left(\rho-\rho',\theta-\theta'\right) \, d\theta' d\rho'.
\end{align*}
It is shown in \cite{andersson2005fast} that $\Rad^\# f(x) = 2\Radlp^\# f(\theta,\rho)$, where the factor 2 comes from the fact that the corresponding integration in the polar representation $(\theta,s)$ is only done in the half-plane $s>0$. The log-polar back-projection operator has the same problem as the log-polar Radon transform in dealing with $s=0$, and in a similar fashion we make use of \emph{partial back-projections} in order to avoid this problem. Because of the relation \eqref{Radon_all_pb} and the fact that the filter operator $\mathcal{W}$ can be applied to each partial Radon data individually, it will be enough to consider partial back-projections for Radon data with $\theta\in\left[-\betah+m\beta,\betah+m\beta\right]$ according to the setup of Figure \ref{fig:spans}. By applying $T_m^{-1}$ to each of the partial back-projection, and summing up the results, we will recover the original function. For detailed calculations, we refer to \cite{andersson2005fast}.

The idea is thus to split the Radon data into $M$ parts, where each part is transformed according to Figure \ref{fig:spans}a). For each part, the filtered data is back-projected according to
\begin{equation}\label{partial_bck_1}
\int_{-\infty}^{\infty} \int_{-\betah}^\betah g(\theta',\rho') \zeta^\# \left(\rho-\rho',\theta-\theta'\right) \, d\theta' d\rho'.
\end{equation}
Since we assumed that our original function had support inside the unit circle, we are only be interested in the contributions inside the disc $D$. Since only lines with $\rho \in [\log(a_r),0]$ go through this circle, the integration in the $\rho$ variable above can be limited to $\rho \in [\log(a_r),0]$. Similarly as for the Radon transform, we now want to write this (finite) convolution as a periodic convolution. Figure  \ref{fig:conv_support}b) illustrates how this can be achieved. The black solid lines show translations of the curves $\rho=-\log(\cos(\theta))$ representing the back-projection integral in the log-polar coordinates. The black dots now show the perimeter of a support of the Radon data, indicated by dark gray in the left illustration of Figure \ref{fig:conv_support}. The dark gray curves the back-projection illustration now show the translations of $\rho=-\log(\cos(\theta))$ that will give a contribution inside the disc $D$. The light gray curves illustrate contributions that fall outside the support of $D$. The red curves show contributions that will fall outside the range $[-\beta,\beta]\times [\log(a_r),0]$, and the blue curves show the effect when these lines are wrapped back in to the domain $[-\beta,\beta]\times [\log(a_r),0]$. We note that the blue curves do not intersect the circle. Hence, we define the partial log-polar  back-projection operator as
\begin{equation}\label{backproj_lpp}
\Radlppart^\# g(\theta,\rho) = \int_{\log(a_r)}^{0} \int_{-\betah}^\betah g(\theta',\rho') \zeta_{\mathrm{per}}^\# \left(\rho-\rho',\theta-\theta'\right) \, d\theta' d\rho'.
\end{equation}
where $\zeta^\#_{\mathrm{per}}$ is the periodic extension of $\zeta^\#$ defined on $[\log(a_r),0]\times [-\beta,\beta]$, and note that for  $(\theta,\rho)$ corresponding to points inside the domain $D$, it holds that
\begin{align*}
\Radlppart^\# g(\theta,\rho) = \Radlp^\# g(\theta,\rho) .
\end{align*}
We then have that
$$
2\sum_{m=0}^{M-1}  T_m^{-1} \Radlppart^\# \mathcal{\Jlp} \Radlppart T_m f(x)= f(x)
$$
for all $x$ in the unit disc.

\section{Fast evaluation of the log-polar Radon transform and the log-polar back-projection}

Let $h$ be a continuous function on some rectangle $R$ in $\mathbb{R}^2$, and let $\mu$ be a finite measure on $R$, and let $\mu_{\mathrm{per}}$ be its periodic extension. Along the same lines as \cite[Theorem 11.6-3]{zemanian}, it is easy to see that the Fourier coefficients of their periodic convolution satisfies $$\widehat{h*\mu_{\mathrm{per}}}=\hat{h}\hat{\mu},$$
where $\hat{h}\hat{\mu}$ is the pointwise multiplication of the respective Fourier coefficients with respect to $R$. We will use this formula and FFT to fast evaluate $\Radlppart$ \eqref{lograd_periodic} and $\Radlppart^\#$ \eqref{backproj_lpp}.
We use the notations
\begin{align}
f(\theta,\rho) &= \sum_{k_\theta,k_\rho} \widehat{f}_{k_\theta,k_\rho} e^{2\pi i \left( \frac{\theta k_\theta}{2\beta} + \frac{\rho k_\rho}{-\log{a_r}}\right) }, \label{fFTrep}\\
g(\theta,\rho) &= \sum_{k_\theta,k_\rho} \widehat{g}_{k_\theta,k_\rho} e^{2\pi i \left( \frac{\theta  k_\theta}{2\beta} + \frac{\rho k_\rho}{-\log{a_r}}\right) },\label{gFTrep}\\
\zeta(\theta,\rho) &= \sum_{k_\theta,k_\rho} \widehat{\zeta}_{k_\theta,k_\rho} e^{2\pi i \left( \frac{\theta k_\theta}{2\beta} + \frac{\rho k_\rho}{-\log{a_r}}\right) },\label{zetaFTrep}\\
\zeta^\#(\theta,\rho) &= \sum_{k_\theta,k_\rho} \widehat{\zeta^\#}_{k_\theta,k_\rho} e^{2\pi i \left( \frac{\theta k_\theta}{2\beta} + \frac{\rho k_\rho}{-\log{a_r}}\right) }.\label{zetabFTrep}
\end{align}
The Fourier coefficients for the two distributions $\zeta$ and $\zeta^\#$ are given by
\begin{align}
\zeta_{k_\theta,k_\rho}&=\int_{-\beta}^{\beta}\int_{\log(a_r)}^{0}\delta(\cos(\theta)-e^{\rho})e^{-2\pi i\frac{\theta k_\theta}{2\beta}}e^{-2\pi i \frac{\rho k_\rho}{-\log(a_r)}}d\rho d\theta=\nonumber\\
&=\int_{-\beta}^{\beta}\int_{\log(a_r)}^{0}\delta(\cos(\theta)-e^{\rho})e^{-2\pi i\frac{\theta k_\theta}{2\beta}} \left( e^{\rho} \right)^{-2\pi i \frac{k_\rho}{-\log(a_r)}}d\rho d\theta=\nonumber\\
&=\label{exprcoef1}\int_{-\beta}^{\beta}e^{-2\pi i\frac{\theta k_\theta}{2\beta}}(\cos(\theta))^{-2\pi i \frac{ k_\rho}{-\log(a_r)}-1}d\theta,\\
\zeta_{k_\theta,k_\rho}^{\#}&=\int_{-\beta}^{\beta}\int_{0}^{-\log(a_r)}\delta(e^{\rho}\cos(\theta) - 1)e^{-2\pi i\frac{\theta k_\theta}{2\beta}}e^{-2\pi i \frac{\rho k_\rho}{-\log(a_r)}}d\rho d\theta=\nonumber\\
&=\int_{-\beta}^{\beta}\int_{0}^{-\log(a_r)}\delta(e^{\rho}\cos(\theta) - 1)e^{-2\pi i\frac{\theta k_\theta}{2\beta}} \left( e^{\rho} \right)^{-2\pi i \frac{\rho k_\rho}{-\log(a_r)}}d\rho d\theta=\nonumber\\
&=\label{exprcoef2}\int_{-\beta}^{\beta}e^{-2\pi i\frac{\theta k_\theta}{2\beta}}(\cos(\theta))^{-2\pi i \frac{ k_\rho}{\log(a_r)}}d\theta.
\end{align}
Both the integrals on the right hand sides of \eqref{exprcoef1} and \eqref{exprcoef2} are of the form
\begin{align}\label{Pfor}
P(\mu,\alpha,\beta)=\int_{-\beta}^{\beta}e^{i\mu\theta}\cos(\theta)^\alpha d\theta,
\end{align}
where $\mu={-2\pi \frac{k_\theta}{2\beta}}$, $\alpha={-2\pi i \frac{ k_\rho}{-\log(a_r)}-1}$ for \eqref{exprcoef1}  and $\alpha={-2\pi i \frac{ k_\rho}{\log(a_r)}}$ for \eqref{exprcoef2}, respectively. It turns out that there is a closed form expression for the integral \eqref{Pfor}, namely
\begin{align}\label{P_closed_form}
&P(\mu,\alpha,\beta)=
\frac{\Gamma(\frac{\alpha+1}{2})\Gamma(\frac{1}{2})\Gamma(\frac{\alpha+2}{2})}{\Gamma(\frac{\alpha+\mu}{2}+1)\Gamma(\frac{\alpha-\mu}{2}+1)}+\\
&\frac{2\mu\cos(\beta)^{\alpha+2}\sin(\mu\beta)}{(\alpha+1)(\alpha+2)}\,_3F_2\left(1,\frac{\alpha}{2}+\frac{\mu}{2}+1,\frac{\alpha}{2}-\frac{\mu}{2}+1;\frac{\alpha+3}{2},\frac{\alpha}{2}+2;\cos(\beta)^2\right)-\\
&\frac{2\cos(\beta)^{\alpha+1}\cos(\beta \mu)\sin(\beta)}{(\alpha+1)} \,_3F_2\left(1,\frac{\alpha}{2}+\frac{\mu}{2}+1,\frac{\alpha}{2}-\frac{\mu}{2}+1;\frac{\alpha+3}{2},\frac{\alpha}{2}+1;\cos(\beta)^2\right).
\end{align}
In the appendix we derive this expression. However, as it includes gamma and hypergeometric functions with complex arguments, special care has to be taken when evaluating \eqref{P_closed_form} numerically, as cancelation effect easily can cause large numerical errors. For the case where there are not suitable numerical routines readily available for evaluation of these special functions, we briefly describe an alternative way to evaluate \eqref{Pfor} numerically. Note that for $\mu={-2\pi \frac{k_\theta}{2\beta}}$, the integral is well suited for evaluation by FFT, since for a fixed value of $k_\rho$ (fixed $\alpha$) we can obtain $\zeta_{k_\theta,k_\rho}$ and $\zeta^\#_{k_\theta,k_\rho}$ for all integers $k_\theta$ in a given range by evaluating the integral of \eqref{Pfor} by the trapezoidal rule by means of FFT. However, for this procedure to be accurate, we need to oversample the integral, and make use of end-point corrections. The function $\cos(\theta)^\alpha$ will be oscillatory, but the oscillation is determined by the fixed parameter $\alpha$. Neglecting boundary effects, we can therefore expect the trapezoidal rule to be efficient, provided that sufficient oversampling is used.

For the boundary effects, we can make use of end-point correction schemes \cite{alpert1995high}. For the computations used in this paper, we have used an oversampling factor of $8$ and an 8-order end-point correction with weights
\begin{align}
1+\frac{1}{120960} \Big[-23681, 55688, -66109, 57024, -31523, 9976, -1375\Big] \label{weightsa}.
\end{align}

Suppose next that we only know values of $f$ and $g$ in \eqref{lograd_periodic} and \eqref{backproj_lpp}, respectively, on an equally spaced sampling covering $[\log(a_r),0]\times [-\beta,\beta]$, i.e., $f$ and $g$ are known on the lattice
\begin{equation}\label{rhothetalattice}
\left\{ \frac{j_\theta}{2\beta N_\theta}, \frac{j_\rho}{-\log(a_r) N_\rho} \right\},\quad -\frac{N_\theta}{2} \le j_\theta < \frac{N_\theta}{2}, \quad 0 \le j_\rho < N_\rho,
\end{equation}
with $N_\theta=2\left\lceil \frac{N_\theta}{2M} \right\rceil$.
We denote these values by $f_{j_\theta,j_\rho}$ and $g_{j_\theta,j_\rho}$, respectively.
In order for \eqref{lograd_periodic} and \eqref{backproj_lpp} to be meaningful, we need to have continuous representations of $f$ and $g$. This is particularly important as  $\zeta$ and $\zeta^\#$ are distributions. A natural way to do this, is to define $\widehat{f}_{k_\theta,k_\rho}$ and $\widehat{g}_{k_\theta,k_\rho}$ by the discrete Fourier transform of $f_{j_\theta,j_\rho}$ and $g_{j_\theta,j_\rho}$, i.e.,  let
\begin{equation} \label{fhatdef}
\widehat{f}_{k_\theta,k_\rho} =
\begin{dcases}
\sum_ { j_\theta =-\frac{N_\theta}{2} }^{\frac{N_\theta}{2}-1} \sum_{j_\rho=0}^{N_\rho-1} f_{j_\theta,j_\rho} e^{-2\pi i \left( \frac{j_\theta k_\theta}{N_\theta} + \frac{j_\rho k_\rho}{N_\rho} \right) },& \mbox{if} \quad -\frac{N_\theta}{2} \le k_\theta < \frac{N_\theta}{2}, \quad 0 \le k_\rho < N_\rho,\\
0 & \mbox{otherwise,}\\
\end{dcases}
\end{equation}
and
\begin{equation} \label{ghatdef}
\widehat{g}_{k_\theta,k_\rho} =
\begin{dcases}
\sum_ { j_\theta =-\frac{N_\theta}{2} }^{\frac{N_\theta}{2}-1} \sum_{j_\rho=0}^{N_\rho-1} g_{j_\theta,j_\rho} e^{-2\pi i \left( \frac{j_\theta k_\theta}{N_\theta} + \frac{j_\rho k_\rho}{N_\rho} \right) },& \mbox{if} \quad -\frac{N_\theta}{2} \le k_\theta < \frac{N_\theta}{2}, \quad 0 \le k_\rho < N_\rho,\\
0 & \mbox{otherwise.}\\
\end{dcases}
\end{equation}
We can then use \eqref{fFTrep} and \eqref{gFTrep} \textit{to define} continuous representations of $f$ and $g$. Using this approach, the values of \eqref{lograd_periodic} and \eqref{backproj_lpp} are also well defined, namely
$$
\Radlppart f(\theta,\rho) = \sum_{k_\theta,k_\rho} \widehat{f}_{k_\theta,k_\rho}  \widehat{\zeta}_{k_\theta,k_\rho} e^{2\pi i \left( \frac{\theta k_\theta}{2\beta} + \frac{\rho k_\rho}{-\log{a_r}}\right) }, \label{Rad_lp_discrete}\\
$$
and
$$
\Radlp^\# g(\theta,\rho) =\sum_{k_\theta,k_\rho} \widehat{g}_{k_\theta,k_\rho} \widehat{\zeta^\#}_{k_\theta,k_\rho} e^{2\pi i \left( \frac{\theta k_\theta}{2\beta} + \frac{\rho k_\rho}{-\log{a_r}}\right) }\label{Rad_bck_lp_discrete}.
$$
In the case where the two transforms are to be evaluated at $(\theta,\rho)$ on the lattice \eqref{rhothetalattice}, the corresponding sums above can be rapidly evaluated by used FFT.
\section{Sampling rates}
There are three different cases for which it is natural to use equally spaced discretization: For the representation of $f$ in Cartesian coordinates $(x_1,x_2)$ covering the unit circle $S^1$;
for the polar representation $(\theta,s)$ of the sinograms $\Rad(f)$; and for the log-polar coordinates $(\theta,\rho)$ for evaluation of $\Radlp$ and $\Radlp^{\#}$. We will use a rectangular grid in all three coordinate systems and interpolate data between them as we work with the different domains.

In this section we derive guidelines for how to choose discretization parameters. These will be based on the assumption that $\hat{f}$ is ``essentially supported'' in a disc with radius $N/2$, in the sense that contributions from the complement of this disc can be ignored without affecting accuracy of our calculations. Due to the uncertainty principle, $\hat {f}$ can not have its support included in this disc, since $f$ itself is supported in a disc with radius $1/2$, but it may work quite well in practice and is therefore still convenient to use for deriving sampling rates.\footnote{More strict results can be obtained by introducing notions such as numerical support. Recall that any compactly supported $L^1$ function has compact numerical support in both $x$ and $\xi$ by the Riemann-Lebesgue lemma. For instance, Gaussians have (small) compact numerical support in both the spatial and the frequency domain. The presentation will quickly become substantially more technical and we therefore follow the signal processing practice and use Nyquist sample rates despite not having infinitely long samples.} We also base our arguments on refinements of the Nyquist sampling rate, more precisely the Paley-Wiener-Levinson theorem.

In the spatial variables $(x_1,x_2)$, the Nyquist sampling rate (corresponding to the assumption on the support of $\hat{f}$) is $1/N$. This leads us to cover the unit circle with a grid of size $N\times N$. We use
\begin{equation}\label{Xgrid}
X= \left\{ \frac{j_1}{N}, \frac{j_2}{N}\right\}, \quad -\frac{N}{2} \le j_1,j_2 < \frac{N}{2}.
\end{equation}

To derive recommendations for the sampling in the polar coordinates $(\theta,s)$ (for the entire $\Rad f$), recall the Fourier slice theorem, which states that
\begin{align*}
&\int_{-\infty}^\infty \Rad f (s,\theta) e^{-2\pi i s \sigma} \, ds = \int_{-\infty}^\infty \int_{-\infty}^\infty f (s\theta+t\theta^\perp) e^{-2\pi i s \sigma} \, ds \, dt =\\
&\int_{-\infty}^\infty \int_{-\infty}^\infty f (x) e^{-2\pi i \sigma (x\cdot \theta) } \, ds \, dt=\int_{-\infty}^\infty \int_{-\infty}^\infty f (x) e^{-2\pi i (\sigma \theta) \cdot x } \, ds \, dt =\widehat{f}(\sigma \theta).
\end{align*}
Thus $\Rad f(\theta,s)=\mathcal{F}^{-1}_{\sigma\rightarrow s}(\hat{f}(\sigma \theta))$, and hence the Nyquist sampling rate for $s$ is $\triangle s=\frac{1}{N}$, yielding that $N_s = N$.
Let $\triangle \thetap$ denote the angular sampling rate in the polar coordinate system. Since $f$ is supported in the unit disc, the Nyquist sampling rate in the frequency domain $\xi$ is 1 (on a rectangular infinite grid). However, the multidimensional version of the Paley-Wiener-Levinson theorem roughly says that it is sufficient to consider an irregular set of sample points whose maximal internal distance between neighbors is 1. Since we have assumed that the values of $\widehat{f}(\xi)$ for $|\xi| >\frac{N}{2}$ are negligible, this leads us to choose $\triangle \thetap $ and $\triangle \sigma$ (the latter will not be used) so that the polar grid-points inside this circle has a maximum distance of 1. It follows that
\begin{equation}\label{delta_theta_org}
\triangle \thetap = \frac{2}{N}.
\end{equation}
Since we want to cover an angle span of $[0,\pi]$, this leads to $N_\theta \approx \frac{\pi}{2} N$. We denote the polar grid by
\begin{equation}\label{Polargrid}
\Sigma= \left\{ j_\theta \triangle \thetap,  j_s \triangle s \right\}, \quad 0 \le j_\theta < N_\theta, \quad -\frac{N}{2} \le j_s < \frac{N}{2}.
\end{equation}

    Practical tomographic measurements can therefore typically have a ratio $\frac{N_\theta}{N_s} \approx 1.5$ between the sampling rates in the $\theta$- and $s$-variables. We refer to \cite{kak1988principles, natterer1986computerized} for more details, and proceed to discuss the sampling in the log-polar coordinates.

\begin{figure}
\centering
\subfloat[][]{\includegraphics[trim = 0mm 0mm 0mm 0mm,clip=true,width=0.47\textwidth]{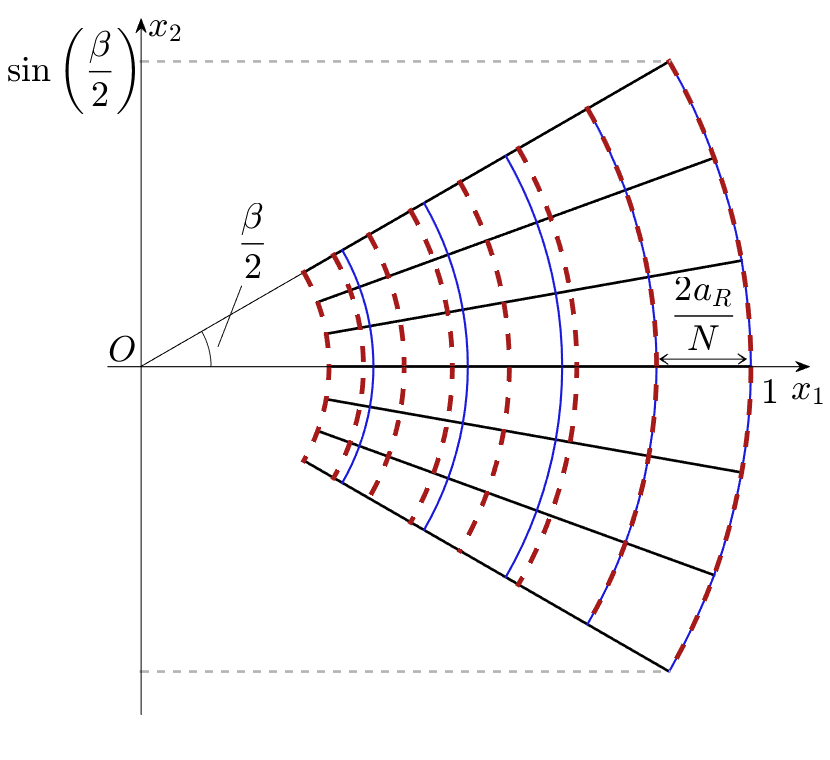}}\hspace{2mm}
\subfloat[][]{\includegraphics[trim = 34mm 20mm 4mm 0mm,clip=true,width=0.5\textwidth]{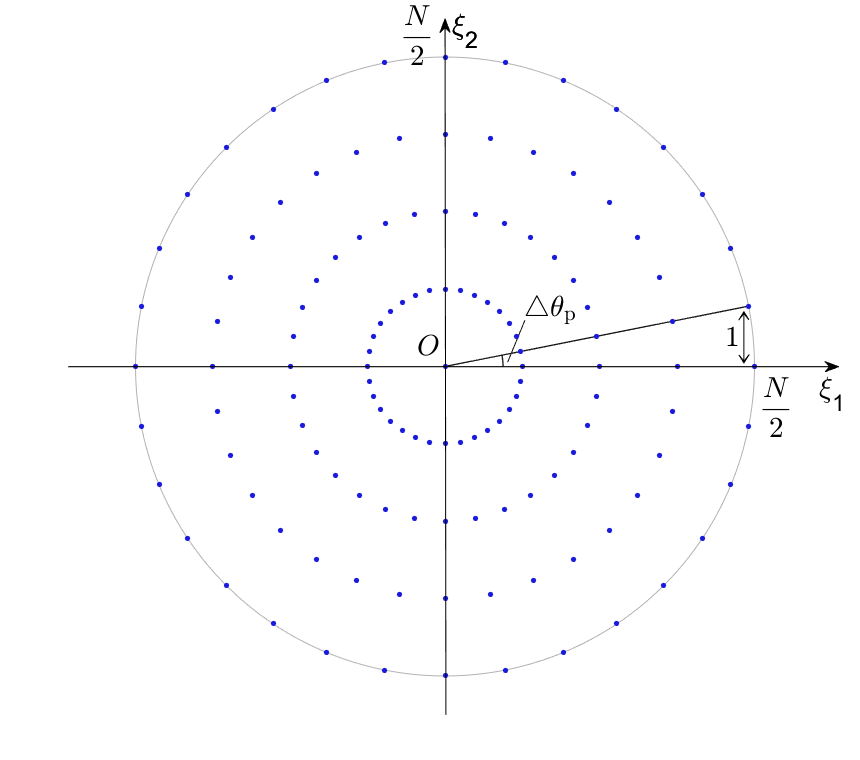}}
\caption{(a) Log-polar grid. Samples in the $\rho$ variable are chosen in order not to lose the accuracy of measuring data. (b) The bandwidth of a partial back-projection function. }
\label{fig:gridslp}
\end{figure}

To distinguish between the coordinate sampling parameters we use $\triangle \thetalp$ for angular sampling rate in the log-polar coordinates. Recall that we wish to accurately represent functions of the form $T_m f$, which is a rotation, dilation (with factor $2a_R$) and translation of $f$. The support of $T_m f$ lies in the grey circle of Figure \ref{fig:spans} and its essential frequency support is then inside the disc of radius $\frac{N}{4a_R}$. 

Figure \ref{fig:gridslp} illustrates the setup where the black lines indicate equally spaced samples of $\thetalp$; the blue curves indicate equally spaced sampling in $\slp$; and where the red dashed curves indicate equally spaced sampling in $\rho$. Note that the maximum distance between points in the $x_1$ direction occurs when $s=1$ or $\rho=0$. It is thus clear that $\triangle \slp$ should equal $\triangle x_1=\frac{2 a_R}{N}$, whereas $\triangle \rho$ is determined by
\begin{align*}
\max_{\rho \in [\log(a_r),0]} e^{\rho}-e^{\rho-\triangle \rho} = \max_{\rho \in [\log(a_r),0]} e^{\rho}\left(1- e^{-\triangle \rho}\right) \le\triangle x_1.
\end{align*}
As the largest distance occurs when $\rho=0$ it follows that
$$
\triangle \rho \le -\log\left(1- \frac{2a_R}{N}\right),
$$
and consequently, since the total distance that is to be covered is $-\log(a_r)$, that
\begin{align}\label{vrho}
&N_\rho\ge\left\lceil\frac{\log(a_r)}{\log(1-\frac{2a_R}{N})}\right\rceil.
\end{align}
where the notation $\lceil x \rceil$ denotes the nearest integer greater than or equal to $x$.

For the determination of sample rate in the angular variable for the representation of $T_m f$ in the disc $D$, we have that
$$
\triangle \thetalp\approx \sin(\triangle \thetalp) \le \frac{2 a_R}{N}.
$$
Hence, we introduce
\begin{equation}\label{Lopolargrid_f}
\Omega_{\mathrm{lp}}= \left\{ j_\theta \triangle \thetalp , j_\rho  \triangle \rho \right\}, \quad -\frac{\beta}{2 \triangle \thetalp} \le j_\theta < \frac{\beta}{2 \triangle \thetalp}, \quad -N_\rho < j_\rho \le 0,
\end{equation}
for the representation of $T_m f$ in $D$. The representation of the partial Radon transform of $T_m f$ needs a reduced sample rate compared to \eqref{Lopolargrid_f}. As the Radon transform is applied as a convolution on a log-polar grid (i.e. by a Fourier multiplier), the higher frequencies in the $\theta$-direction will not be needed.  Hence, we can apply a low-pass filter in the $\theta$-direction, and apply the FFT operation the grid
\begin{equation}\label{Lopolargrid_polar}
\Omega_\mathrm{p}= \left\{ j_\theta \triangle \thetap, j_\rho \triangle \rho \right\}, \quad -\left\lfloor \frac{N_\theta}{2 M} \right\rfloor  \le j_\theta < \left\lfloor \frac{N_\theta}{2 M} \right\rfloor, \quad N_\rho < j_\rho \le 0,
\end{equation}
for the computation of the partial Radon tranform of $T_m f$. The grid $\Omega_\mathrm{p}$ will also be useful when resampling Radon data in log-polar coordinates.

\section{Interpolation}
As mentioned in the previous section it will be natural to use equally spaced sampling in the Cartesian, polar and log-polar coordinate systems, respectively. We typically want to reconstruct data on an equally spaced Cartesian grid; the tomographic data is sampled in equally spaced polar coordinates; and both the Radon transform and the back-projection can be rapidly evaluated by FFT when sampled on an equally spaced log-polar grid. There are several ways to interpolate between these coordinate systems. In this work we will make use of cubic (cardinal) B-spline interpolation. We will discuss how to incorporate some of the interpolation steps in the FFT operations, and also discuss how the B-spline interpolation can be efficiently implemented on GPUs.

The cubic cardinal B-spline is defined as
\begin{align*}
B(x)= \begin{cases}
(3|x|^3-6|x|^2+4)/6, & 0\le|x|<1,\\
(-|x|^3+6|x|^2-12x+8)/6, & 1\le|x|<2,\\
0, & |x| \ge 2.
\end{cases}
\end{align*}
This function is designed so that it is non-negative, piecewise smooth, and $C^2$ at $|x|=0,1,2$. Note that at integer points, it holds that
\begin{equation}\label{Bcoeff}
B(j) = \begin{cases}
\frac{2}{3} & \mbox{if } j=0, \\
\frac{1}{6} & \mbox{if } |j|=1, \\
0 & \mbox{otherwise.}
\end{cases}
\end{equation}
When used as a filter, it will smooth out information and is thus acting as a low-pass filter. Consequently, it can not be used directly for interpolation.

The Fourier series given by the coefficients in \eqref{Bcoeff}, which we denote by $\widehat{B}$, is given by
\begin{equation}\label{Bhat}
\widehat{B}(\xi) = \sum_{j} B(j) e^{-2\pi i j \xi} = \frac{1}{6} \left( e^{2\pi i x \xi} + 4 +  e^{-2\pi i x \xi}\right) = \frac{2}{3} + \frac{1}{3} \cos(2\pi \xi).
\end{equation}


Suppose that equally spaced samples $f_k$ of a function $f$ (in one variable) are available. We want to recover values of $f$ at arbitrary points $x$ using
\begin{equation}\label{Bsplineinterp}
f(x)=\sum_{k} (Q f)_k B\left(x-\frac{k}{N}\right).
\end{equation}
Since $B$ has short support, only values of $(Q f)_k$ for $k\approx x N$ will contribute in this sum. The operator $Q$ is a pre-filter operation, that is compensating for the fact that convolution with $B$ suppresses high frequencies. The pre-filter operation is boosting high frequencies in the samples $f_k$. It can be computed in different ways. Perhaps the most direct way is to define $Q$ in the Fourier domain (by the discrete Fourier transform), where it essentially becomes division by $\widehat{B}$ (upon scaling and sampling). In this case, it is easy to see that the convolution with $B$ and the pre-filter operation will cancel each other at points $x=\frac{j}{N}$, i.e., the original function is recovered at the sample points; which is a requisite for any interpolation scheme.

As we will compute the Radon transform and the back-projection by means of FFT in log-polar coordinates, we can in some steps incorporate the pre-filter step $1/\widehat{B}$ in the Fourier domain, at virtually no additional cost. However, not all of the pre-filter operations can be incorporated in this way. While the pre-filter easily can be applied by separate FFT operations, we want to limit the total number of FFT operations, as these will be the most time-consuming part in the implementations we propose.

As an alternative to applying the pre-filter in the Fourier domain, it can be applied by recursive filters. In \cite{ruijters2010gpu} these operations are derived by using the $Z$-transform.

%
%
%
%
It turns out that if we define
$$
(Qf^+)_k=6 f(k)+(\sqrt{3}-2) (Qf^+)_{k-1},
$$
then
\begin{equation}\label{Qfrecursive}
(Qf)_k=(\sqrt{3}-2)((Qf)_{k+1}-(Qf^+)_k),
\end{equation}
cf. \cite[equations (12,13)]{ruijters2010gpu}, where the boundary condition on $Qf$ and $Qf^+$ are given in  \cite[equations (14,15)]{ruijters2010gpu}. These filters can be efficiently implemented on GPUs.

Two-dimensional prefiltration can be done in two steps, one in each dimension. We will use the notation $Q f$ for the pre-filtering also in this case. On a GPU this implies doing operations on rows and columns separately.  However, there are highly optimized routines for transposing data, which means that the prefiltration can be made to act only on column data. This is done to improve the so-called memory coalescing and GPU cache performance \cite{wilt2013cuda}. Memory coalescing refers to combining multiple memory accesses into a single transaction. In this way the GPU threads run simultaneously and substantially increased cache hit ratios are obtained.

Let us now turn our focus to the convolution step in \eqref{Bsplineinterp}. Let
\begin{align}\label{floorx}
&k=\lfloor N x \rfloor,\\
&\nonumber\alpha=x-\lfloor x \rfloor.
\end{align}
The sum \eqref{Bsplineinterp} then reduces to
\begin{align}
&f(x)=w_0(\alpha) (Q f)_{k-1}+w_1(\alpha)(Q f)_{k}+w_2(\alpha)(Q f)_{k+1}+w_3(\alpha)(Q f)_{k+2},\text{ where }\label{wsumn}\\
&w_0(\alpha)=B(\alpha-2),\quad w_1(\alpha)=B(\alpha-1),\quad w_2(\alpha)=B(\alpha),\quad w_3(\alpha)=B(\alpha+1).\nonumber
\end{align}
We now discuss how this sum can be evaluated using linear interpolators. As mentioned previously, linear interpolation is executed fast on GPUs. In \cite{sigg2005fast} it is shown how this can be utilized to conduct efficient cubic interpolation. The cubic interpolation is expressed as two weighted linear interpolations, instead of four weighted nearest neighbor look-ups, yielding $2^d$ operations instead of $4^d$ for conducting cubic interpolation in $d$ dimensions.

We briefly recapitulate the approach taken in \cite{ruijters2008efficient,sigg2005fast}. Given coefficients $(Qf)_k$, let $ Qf_{\mathrm{lin}}$ be the linear interpolator
\begin{align*}
Qf_{\mathrm{lin}}(x)  &=(1-(\alpha))(Qf)_{k}+\alpha (Qf)_{k+1}\\
&=(1-(x-\lfloor x \rfloor))(Qf)_{\lfloor N x \rfloor}+(x-\lfloor x \rfloor) (Qf)_{\lfloor N x \rfloor+1}.
\end{align*}
The sum \eqref{wsumn} can be then be  written as
\begin{align}
f(x)&=(w_0(\alpha)+w_1(\alpha))  Qf_{\mathrm{lin}}\left(k-1+\frac{w_1(\alpha)}{w_0(\alpha)+w_1(\alpha)}\right) \nonumber\\
&+(w_2(\alpha)+w_3(\alpha)) Qf_{\mathrm{lin}}\left(k+1+\frac{w_3(\alpha)}{w_2(\alpha)+w_3(\alpha)}\right). \label{cubic_interp_reduced}
\end{align}
The evaluation of the function $Qf_{\mathrm{lin}}(x) $ can be performed by hard-wired linear interpolation on GPU. In modern GPU architecture the so-called texture memory (cached on a chip) provide effective bandwidth by reducing memory requests to the off-chip DRAM. The two most useful features of these kind of memory with regards to conducting B-spline interpolation, are
\begin{enumerate}
\item The texture cache is optimized for the 2D spatial locality, giving best performance to GPU threads that read texture addresses that are close together.
\item Linear interpolation of neighboring values can be performed directly in the GPUs texture hardware, meaning that the cost for computing the interpolation is the same as reading data directly from memory.
\end{enumerate}
This implies that the cost for memory access in \eqref{cubic_interp_reduced} will be two instead of four as only two function calls of $Qf_{\mathrm{lin}}$ are made in \eqref{cubic_interp_reduced}, and in two dimensions the corresponding reduction from 16 memory access operations to 4 will give significant improvement in computational speed.

\section{Algorithms}
We now have the necessary ingredients to present detailed descriptions on how to rapidly evaluate the Radon transform and the back-projection operator by FFT in log-polar coordinates. In the algorithms below, we let $\widehat{B}_{k_\theta,k\rho}$ denote the values of the two-dimensional counterpart of \eqref{Bhat}, scaled to represent the sampling on $(\theta,\rho)\in\left[-\betah,\betah\right]\times \left[\log a_r,0 \right]$.

\begin{algorithm}
\caption{Fast Radon transform}
\label{Radon_alg}
\begin{algorithmic}[1]
\State Given $f$ sampled at $X$ \eqref{Xgrid} compute $Qf$ by \eqref{Qfrecursive}
\For{$m=0, \dots M-1$}
\State Resample $T_m f$ at $\Omega_\mathrm{lp}$ \eqref{Lopolargrid_f} by \eqref{Bsplineinterp}
\State Downsample from $\Omega_\mathrm{lp}$ \eqref{Lopolargrid_f} to $\Omega_\mathrm{p}$ \eqref{Lopolargrid_polar}
\State Multiply result by $e^\rho$
\State Apply the log-polar Radon transform with pre-filtering incorporated, i.e., compute $\widehat{f}_{k_\theta,k_\rho}$ from \eqref{fhatdef}, and evaluate
$$
\sum_{k_\theta,k_\rho} \widehat{f}_{k_\theta,k_\rho}  \frac{\widehat{\zeta}_{k_\theta,k_\rho}}{\widehat{B}_{k_\theta,k_\rho}} e^{2\pi i \left( j_\theta k_\theta \frac{\triangle \thetap}{2N_\theta} + j_\rho k_\rho \frac{\triangle \rho}{N_\rho}\right) }
$$
by using FFT.
\State Resample from $\mathsf{S}_m^{-1} \Omega_\mathrm{p}$ \eqref{Lopolargrid_polar} to $\Sigma$ \eqref{Polargrid} by \eqref{Bsplineinterp}.
\EndFor
\end{algorithmic}
\end{algorithm}

\begin{algorithm}
\caption{Fast back-projection}
\label{BP_alg}
\begin{algorithmic}[1]
\State Given $g$ sampled at $\Sigma$ \eqref{Polargrid} compute $Qg$ by \eqref{Qfrecursive}
\For{$m=0, \dots M-1$}
\State Resample $g(\mathsf{S}_m)$ at  $\Omega_\mathrm{p}$ \eqref{Lopolargrid_polar} by \eqref{Bsplineinterp}
\State Apply the log-polar back-projection with pre-filtering incorporated, i.e., compute $\widehat{g}_{k_\theta,k_\rho}$ from \eqref{ghatdef} and evaluate
$$
\sum_{k_\theta,k_\rho} \widehat{g}_{k_\theta,k_\rho}  \frac{\widehat{\zeta^\#}_{k_\theta,k_\rho}}{\widehat{B}_{k_\theta,k_\rho}} e^{2\pi i \left( j_\theta k_\theta \frac{\triangle \thetap}{2N_\theta} + j_\rho k_\rho \frac{\triangle \rho}{N_\rho}\right) }
$$
by using FFT.
\State Resample from $\mathsf{T}_m^{-1} \Omega_\mathrm{p}$ \eqref{Lopolargrid_polar} to $X$ \eqref{Xgrid} by \eqref{Bsplineinterp} to obtain partially back-projected data.
\EndFor
\State Sum up the $M$ partial back-projections.
\end{algorithmic}
\end{algorithm}

Let us end with a few remarks on time complexity. The most time consuming part of both algorithms are the convolutions that are implemented by FFT. In total,  $2 M$ FFT operations needs to be computed (including forward and backward FFT's), and each operation will be done on a grid of size $ \frac{2 N_\theta}{M} \times  N_\rho$.  For $M=3$, it follows from \eqref{vrho} that
$$
N_\rho \approx -\frac{N \log(a_r)}{2a_R} \approx 2.1 N,
$$
and $N_\rho$ is monotonically decreasing for increasing $M$. It thus hold approximately that twice as many samples in $\rho$ compared to that originally used for the sampling of $s$ ($N_s=N$), and in the angular variable we also need to sample with twice as many parameters in order to avoid aliasing effects. In total we need to use an oversampling of about 4 in the FFT operations. Note that this is the same oversampling that is generally needed for computing the convolution of two functions if aliasing is to be avoided. The total cost for applying the Radon transform and the back-projection operator is thus the same as would be expected for a generic convolution.

In addition to the FFT operations, interpolation and one-dimensional filter operations are needed in the implementations. In our simulations these operations will typically take about 25\% of the total time.
\section{Performance and accuracy tests}
Let us start by briefly discussing the filters used in the filter back-projection \eqref{FBP}. This discussion is included in order to better interpret the errors obtained when comparing different methods. For theoretically perfect reconstruction (with infinitely dense sampling) the filter $\mathcal{W}$ in \eqref{FBP} is given by
\begin{figure}
\centering
\subfloat{\includegraphics[clip=1,trim = 12mm 2mm 12mm 25mm,clip=true,width=0.5\textwidth]{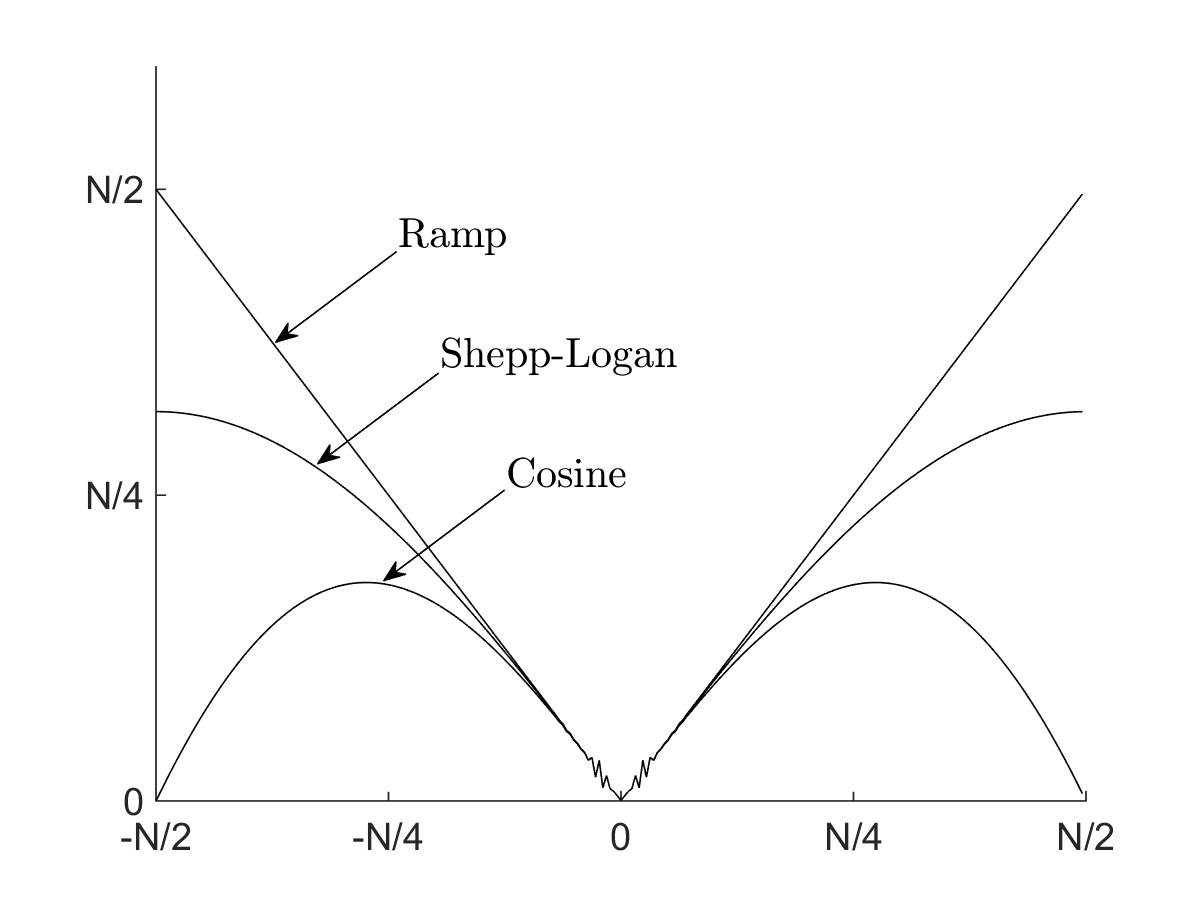}}
\caption{Filters for Radon transform inversion.}
\label{fig:filters}
\end{figure}
\begin{equation} \label{ramp}
\widehat{w}_{\mathrm{ramp}}(\sigma) = |\sigma|.
\end{equation}
This filter is sometimes referred to as the \emph{ramp filter}. If the sampling rate is insufficient in relation to the frequency content of $f$ (or the object upon which the measurements is conducted on), it can be desirable to suppress the highest frequencies in order to localize the effects of the insufficient sampling. There is a relation between one-dimensional filtering in the $s$-direction of Radon data, and a two dimensional convolution in the spatial domain \cite[p.102]{natterer1986computerized}. It can be explained by the Fourier slice theorem,  which describes how Radon data can be converted to a polar sampling in the frequency domain by taking one-dimensional Fourier transforms in the $s$-direction. The application of the one-dimensional filter will then correspond to a change in amplitude (and possible phase) along the lines (indicated by dots with the same angles) in Figure \ref{fig:gridslp}b). As there is no information outside the circle with radius $N/2$, the action of the ramp filter will be equivalent of applying a two-dimensional convolution to the original function $f$ using a two-dimensional filter with Fourier transform
\begin{equation}\label{Framp}
\widehat{W}_{\mathrm{ramp}}(\xi) = \begin{cases}
1 & \mbox{if } |\xi| \le \frac{N}{2}, \\
0 & \mbox{otherwise}.
\end{cases}
\end{equation}
We see that if $f$ contains more high frequent information than the one prescribed by the sampling rate $N$, then the sharp cutoff in \eqref{Framp} can yield artifacts, and in the presence of noise in the Radon sampling the high-frequency boosting of \eqref{ramp} will boost the noise. Replacing the ramp filter \eqref{ramp} with a filter that goes smoothly to zero at the highest (discrete) frequencies, will thus yield an image that is slightly smoother, but on the other hand an image with suppressed high-frequency noise and artifacts due to incomplete sampling. Sometimes, the ramp filter is modified so that it does not reach zero but only reduces the high-frequency amplitudes. Two common choices of filters are the cosine and the Shepp-Logan filters, defined by
\begin{align}
&\label{cosfilt}\widehat{w}_{\cos}(\sigma) = \begin{cases}|\sigma|\cos\left(\frac{2 \pi \sigma}{N}\right),\\
0 & \mbox{otherwise}.
\end{cases} \\
&\widehat{w}_{\mathrm{SL}}(\sigma) = \begin{cases}|\sigma| \sinc\left(\frac{\sigma}{N}\right)\\
0 & \mbox{otherwise}.
\end{cases}. \label{SLfilt}
\end{align}
The cutoff is made above the sampling bandwidth, as it will illustrate the practical effect that the filters have on measured data. The three filters are illustrated in the frequency domain in Figure \ref{fig:filters}. The two-dimensional filters associated with these filters have Fourier representations
\begin{equation*}
\widehat{W}_{\cos}(\xi)=
\begin{cases}
\cos\left(\frac{2 \pi|\xi|}{N}\right) & \mbox{if } |\xi| \le \frac{N}{2}, \\
0 & \mbox{otherwise}.
\end{cases}\label{filters2d}
\end{equation*}
and
\begin{equation*}
\widehat{W}_{\mathrm{SL}}(\xi)=
\begin{cases}
\sinc \left(\frac {|\xi|}{N} \right) & \mbox{if } |\xi| \le \frac{N}{2}, \\
0 & \mbox{otherwise}.
\end{cases}
\end{equation*}
respectively. As the $\widehat{w}_{\cos}$ goes to zero at the highest sampling rate, we can expect smaller error when using this filter compared to the others. On the other hand, the highest frequencies are suppressed, and the results reconstructions will look slightly less sharp.

To illustrate the accuracy of the suggested implementation, we conduct some examples on the \emph{Shepp-Logan} \cite{shepp1974fourier} phantom. We use the modified version introduced in \cite[Appendix B.2]{toft1996radon}.
The function used ($f$) is illustrated in the left panel of Figure \ref{fig:phantom}. The phantom consists of linear combinations of characteristic functions of ellipses, and its support is inscribed in the unit circle. Since the Radon transform is a linear operator, $\Rad f$ is a linear combination of Radon transforms of characteristic functions of ellipses. Since the Radon transform of the characteristic function of a circle can be computed analytically, analytic expressions are available for the Radon transform of the Shepp-Logan phantom by applying transform properties of shifting, scaling and rotation. This Radon transform is depicted in the right panel of figure \ref{fig:phantom}.

The presence of a high-frequency discontinuity caused by the filter can cause artifacts, but also the discontinuity of the derivative at $\sigma=0$. To avoid artifacts from this part, we apply end-point trapezoidal corrections. The effect of this correction can be seen around $\sigma=0$ in Figure \ref{fig:filters}. Our aim is to eliminate as much of the errors as possible to distinguish the errors that the resampling between the different coordinate systems used in the proposed methods introduce.

\begin{figure}
\hspace{0.12\linewidth} Ramp	\hspace{0.19\linewidth} Shepp-Logan \hspace{0.18\linewidth} Cosine
\vspace{2.5mm}

\subfloat{\begin{turn}{90}\,\quad\quad\quad Log-polar \end{turn}}\hspace{0.5mm}
\subfloat{\includegraphics[width=0.296\linewidth]{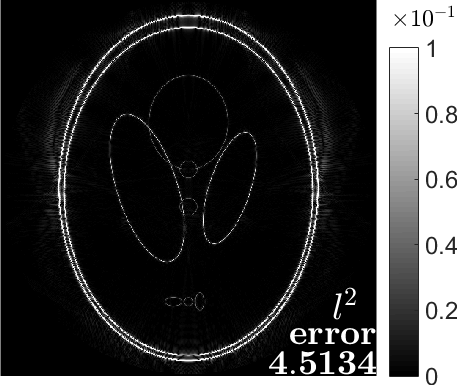}}\hspace{0.5mm}
\subfloat{\includegraphics[width=0.296\linewidth]{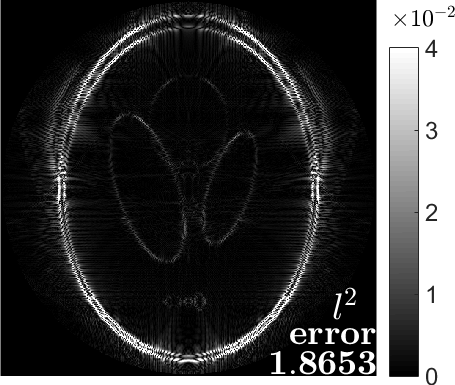}}\hspace{0.5mm}
\subfloat{\includegraphics[width=0.296\linewidth]{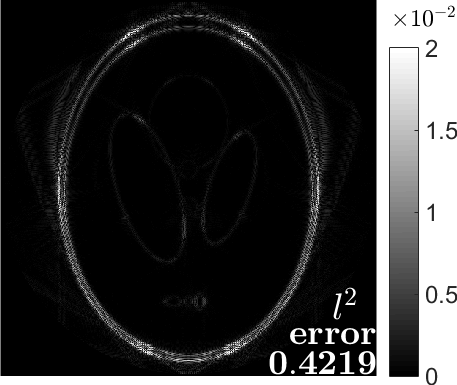}}
\vspace{2.5mm}

\subfloat{\begin{turn}{90}\quad\quad\quad\quad ASTRA\end{turn}}\hspace{1.8mm}
\subfloat{\includegraphics[width=0.296\linewidth]{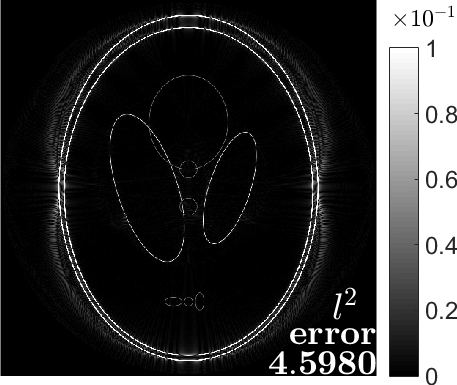}}\hspace{0.5mm}
\subfloat{\includegraphics[width=0.296\linewidth]{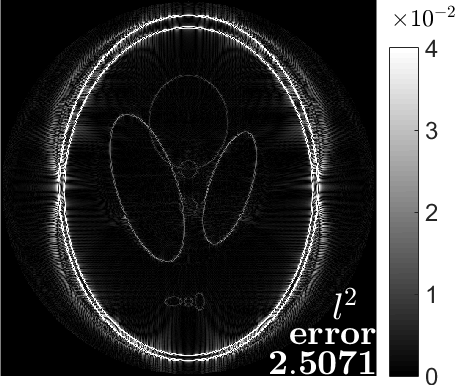}}\hspace{0.5mm}
\subfloat{\includegraphics[width=0.296\linewidth]{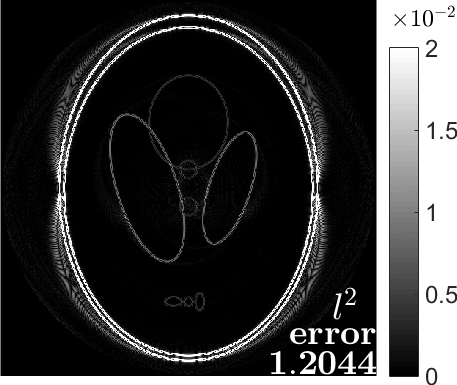}}
\vspace{2.5mm}

\subfloat{\begin{turn}{90}\quad\quad\quad NiftyRec\end{turn}}\hspace{1mm}
\subfloat{\includegraphics[width=0.296\linewidth]{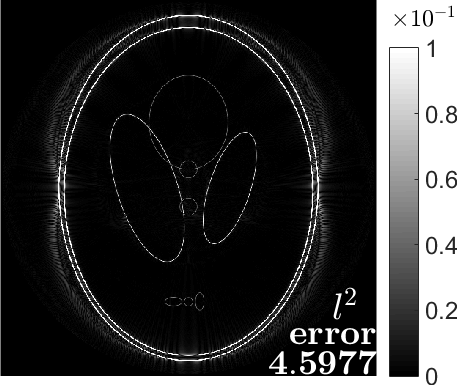}}\hspace{0.5mm}
\subfloat{\includegraphics[width=0.296\linewidth]{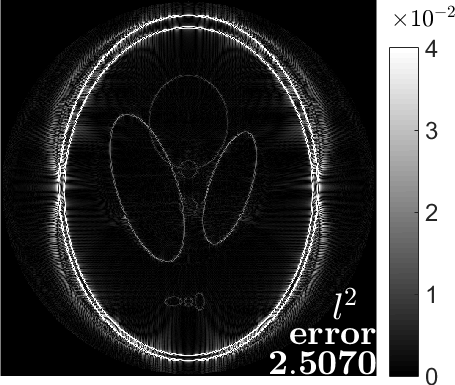}}\hspace{0.5mm}
\subfloat{\includegraphics[width=0.296\linewidth]{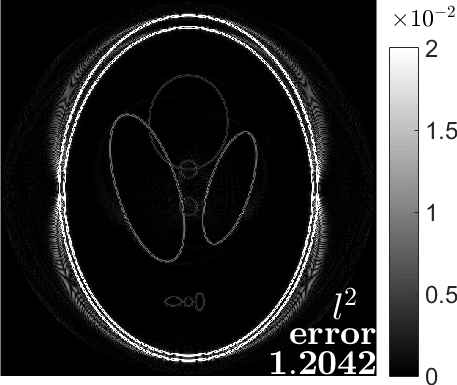}}
\vspace{2.5mm}

\subfloat{\begin{turn}{90}\quad\quad\quad\quad MATLAB\end{turn}}\hspace{1.7mm}
\subfloat{\includegraphics[width=0.296\linewidth]{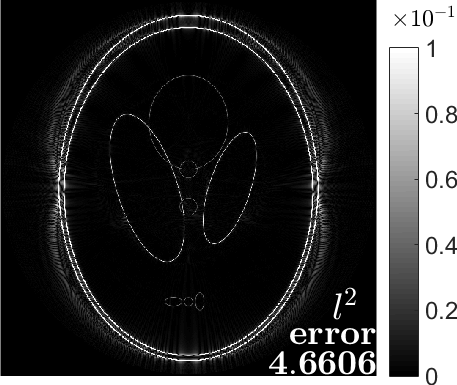}}\hspace{0.5mm}
\subfloat{\includegraphics[width=0.296\linewidth]{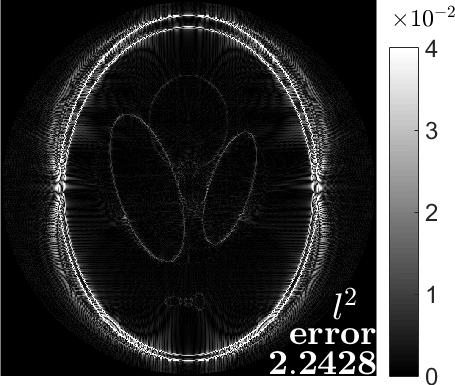}}\hspace{0.5mm}
\subfloat{\includegraphics[width=0.296\linewidth]{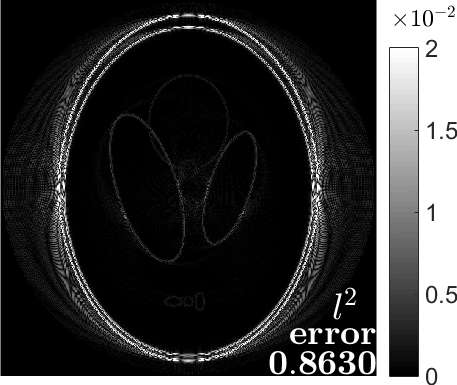}}

\caption{Computational errors of filtered back-projection for the ramp; the Shepp-Logan; and the cosine filters for different methods.}
\label{fig:errors}
\end{figure}

In Figure \ref{fig:errors} we show some reconstruction results from the filtered back-projection using different methods and different filters. For the sake of quality comparison, we use the results from ASTRA Tomography Toolbox \cite{palenstijn2013astra}, NiftyRec Tomography Toolbox \cite{pedemonte2012niftyrec} and the MATLAB image processing toolbox\textsuperscript{TM}. The ASTRA Toolbox uses GPU implementations and the implementation is  described in \cite{palenstijn2011performance,xu2010high}. The NiftyRec Toolbox comes with both GPU and CPU implementations. The toolbox is described in \cite{pedemonte2010gpu,pedemonte20144}.

For the comparisons we have used ASTRA v1.5,  NiftyRec v2.0.1, and MATLAB v2014b.
Both ASTRA and NiftyRec provide MATLAB scripts demonstrating how to use call their routines. For the MATLAB tests, we have used the functions \verb!radon! and \verb!iradon!. Both these functions call compiled routines. To test the back-projection algorithm from Astra toolbox we based our test on the provided routine \verb!s014_FBP!. The iterative tests in the next section are based on the script \verb!s007_3d_reconstruction!. 

The comparison against the NiftyRec toolbox was done by using the back-projection part of the provided demo \verb!tt_demo_mlem_parallel! for parallel beam transmission tomography. The package provides the option to use either GPU or CPU routines. We provide timings for both cases. We use the setup of  \verb!tt_demo_mlem_parallel! for the iterative tests described in the next section.

The reconstructions in the left column of Figure \ref{fig:errors} use the ramp filter \eqref{ramp}; the reconstructions in the middle column use the Shepp-Logan filter \eqref{SLfilt}; and the reconstructions in the right column use the cosine filter \eqref{cosfilt}. The reconstructions were made on a $512\times 512$ grid using 768 samples in the angular direction. For the log-polar reconstruction, we used $M=3$ partial reconstructions, and $N_\rho=2N$. As the cosine filter goes to zero at the boundary, the errors in the reconstruction the right column should better represent the errors caused by the sampling parameters. The reconstruct panels have inscribed $\ell^2$-errors against the filtered versions of the phantom. We can see that ASTRA and NiftyRec produce almost identical errors. We note that the proposed method seems to give the smallest error out of the four methods. The NiftyRec reconstructions refer to the GPU implementation. 

\begin{table}
 \caption{Computational time (in seconds) of the back-projection for sizes $(N_\theta \times N_s)=(\frac{3}{2}N\times N)$ (excluding initialization time).}
 \centering
 \label{table:time_back}
 \begin{tabular}{ | c | c | c | c | c| c | c|}
  \hline
  $N$  & Log-polar & ASTRA & NiftyRec (GPU) & NiftyRec (CPU) & MATLAB \\
  \hline
  256 & 1.6e-03  & 1.0e-02 & 7.1e-02  & 1.3e00 & 3.1e-01  \\
  512 & 6.1e-03  & 3.5e-02 & 5.6e-01  & 1.3e01 & 2.4e00  \\
  1024 & 2.5e-02  & 1.8e-01 & 4.4e00  & 1.2e02 & 1.9e01  \\
  2048 & 9.9e-02   & 1.2e00  & 3.5e+01  & 1.0e03 & 1.5e02  \\
  \hline

 \end{tabular}
\end{table}
In Table \ref{table:time_back} we show times for computing of one back-projection by the different methods mentioned above. To exclude initialization effects and assure high GPU load, the times in Table \ref{table:time_back} are obtained by executing batches of reconstructions (to fill up the GPU memory) in one function call. 

We note that the proposed method is about 5-10 times faster than ASTRA, substantially faster than both the GPU and the CPU implementations of NiftyRec, and up to 1500 times faster than the routines available in the MATLAB toolbox.

For the tests, we have use a standard desktop computer with a Intel Core i7-3820 processor and a NVIDIA GeForce GTX 770 graphic card. We have used NVIDIA CUDA cuFFT library for the FFT operations, and we have used version 7 of the CUDA Toolkit \cite{documentation2015v70}. The computations were performed in single precision. The written GPU program is optimized by using NVidia Nsight profiler, and GPU memory usage, kernels occupancy, instructions fetch and other performance counters were analyzed and improved by using strategies described in \cite{kirk2012programming,sanders2010cuda}.

\section{Iterative methods for tomographic reconstruction}\label{itmethh}
In some situations it is preferable to use an iterative method for doing reconstruction from tomographic data. This could for instance be because of missing data, e.g., that data for some angles are missing; suppression of artifacts \cite{Miqueles:pp5053}; or that additional information about the noise contamination can be used to improve the reconstruction results compared to direction filtered back-projections \eqref{FBP}. Iterative reconstruction methods rely on applying the forward and back-projection operators several times. Iterative algorithms can be computationally expensive when a large number of iterations are required for the algorithm to converge. For that reason, fast algorithms for computing the Radon transform and the associated back-projection can play an important role.

Iterative algebraic reconstruction techniques (ART) are popular tools for reconstruction from incomplete measurements \cite{gubareni2009algebraic}. It aims at solving the set of linear equations determined by the projection data. Transmission based tomographic measurements measure the absorption of a media along a line. This puts a sign constraint on the function we wish to recover. It is clearly not ideal to assume the data is contaminated by normally distributed noise (for which the best estimate is given by the least squares estimates). A more reasonable assumption is that the added noise has a Poisson distribution \cite{barrett1994noise,yan2011expectation}. The simplest iterative method for solving the estimation/reconstruction problem under this noise assumption is by the Expectation-Maximization (EM) algorithm.  The EM-algorithm is well-suited to reconstruct tomography data with non-Gaussian noise \cite{dempster1977maximum,miqueles2014generalized,shepp1982maximum}.  Alternative techniques include for instance the Row Action Maximum Likelihood Algorithm (RAMLA).

More details about the EM-algorithm can be found for instance in \cite{champley2004spect,yan2011expectation}. In our notation, it can (given tomographic data $g$) be expressed as the iterative computations of
\begin{align*}
f^{k+1}=f^{k}\frac{\Rad^\#\left(\frac{g}{\Rad f^k}\right)}{\Rad^\# \chi_C},
\end{align*}
where the function $\chi_C(\theta,s)$ is one if the line parameterized by $(\theta,s)$ passes through the unit disc (the support of $f$) and zero otherwise.
In each step, a Radon transform $\Rad f^k$ and a back-projection $\Rad^\#\left(\frac{g}{\Rad f^k}\right)$ is computed.

It is well-known that a crucial part of GPU computations is host-device memory transfers. For an iterative method such as the one described above, it is possible to keep all the necessary data in the GPU memory, and thereby limit the data copying between host and device memory to an initial guess $f^0$, the measured data $g$, and the final result. As most methods, the proposed method require some initialization steps (e.g., geometry parameters and convolution kernels $\zeta,\zeta^\#$)

\begin{figure}
\centering
\includegraphics[width=0.312\linewidth]{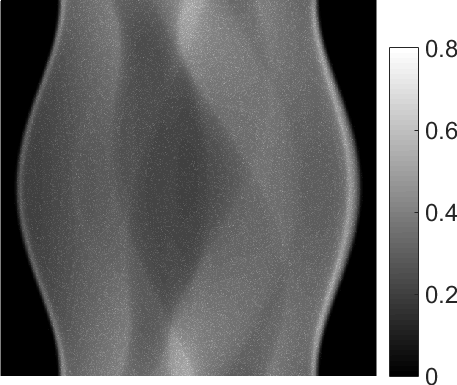}
\includegraphics[width=0.312\linewidth]{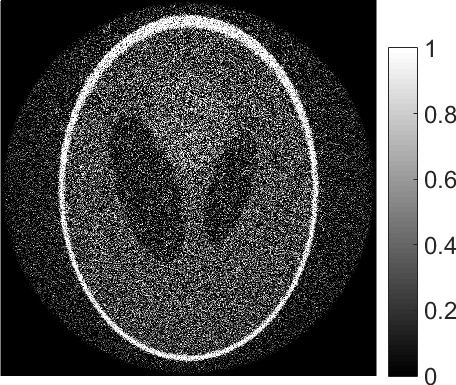}
\includegraphics[width=0.312\linewidth]{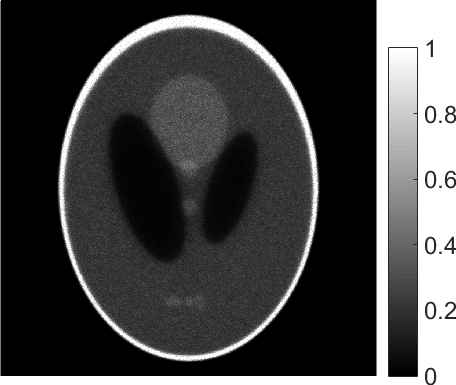}
\caption{Reconstruction of Radon data with noise. Radon data with Poisson noise is shown in the left panel, the result of applying filtered back-projection is shown in the middle panel, the right panel shows the result after 50 iterations of the EM-algorithm.}
\label{fig:EM}
\end{figure}

The obtained GPU program was tested on Radon data with Poisson noise, cf.  Figure \ref{fig:EM}a).  Figure \ref{fig:EM}b) shows the result of applying the filtered back-projection formula, whereas figure \ref{fig:EM}c demonstrates de-noised recovered data after 50 iterations of the EM-algorithm.
In Table \ref{table:timeEM} we present performance results for the EM-algorithm for the proposed method and the methods mentioned in the previous section. The computational times for conducting 100 iterations of the EM-algorithm for different resolution parameters are presented. The Radon data was sampled with parameter $N_\theta = \frac{3}{2} N$, and $N_s=N$. We used $M=3$ partial back-projection for the proposed log-polar based algorithms. The speedup from the previous section is confirmed also for this case.
\begin{table}
\centering
 \caption{Time (in seconds) for 100 iterations of the EM-algorithm for reconstruction 3d data of sizes $(N\times N \times N)$.}
 \label{table:timeEM}
 \begin{tabular}{ | c | c | c | c | c|}
  \hline
$N$& Log-polar & ASTRA & NiftyRec (GPU) & MATLAB \\
  \hline
  256 & 8.8e+01 & 3.1e+02 & 3.6e+03 & 1.6e+04 \\
  512 & 6.9e+02 & 2.3e+03 & 5.8e+04 & 2.6e+05 (*) \\
  1024 & 5.4e+03 & 2.5e+04 & 9.3e+05 & 4.2e+06 (*)\\
  2048 &  4.2e+04 & 3.5e+05 & - & - \\
  \hline
 \end{tabular}

* - estimated by using a reduced number of slices.
\end{table}

\section{Conclusions}
We have described how to implement the Radon transform and the back-projection as convolution operators in log-polar coordinates efficiently on GPUs. We present sampling conditions; provide formulas and numerical guidelines for how to compute the kernels associated with the Radon transform and the back-projection operator; and discuss how the convolutions can be rapidly evaluated by using FFT. The procedure involves several steps of interpolation between data in the Radon domain; the spatial domain; and the log-polar domain. It is comparatively favorable to conduct interpolation in these domains compared to conducting interpolation in the Fourier domain, as the functions will tend to be less oscillatory there. We use cubic spline interpolation which can be efficiently implemented on GPUs, and optimized routines for FFT on the GPU. We conduct numerical tests and see that that we obtain at least as accurate results as that produced by other software packages, but obtained at a substantially lower computational cost.

\section*{Acknowledgments} 
This work was supported by the Crafoord Foundation (20140633) and the Swedish Research Council (2011-5589).

\section*{Appendix} \label{app:Fzeta}
In this appendix, we present the exact result of the integration in \eqref{Pfor}. By applying two partial integrations, it follows that
\begin{align*}
&-\mu^2 P(\mu,\alpha,\beta)=\\&
\int_{-\beta}^{\beta}\left(\frac{d^2}{d\varphi^2}e^{i\mu\varphi}\right)\cos(\varphi)^\alpha d\varphi=\left(\frac{d}{d\varphi}e^{i\mu\varphi}\right)\cos(\varphi)^\alpha |_{-\beta}^{\beta}-\int_{-\beta}^{\beta}\left(\frac{d}{d\varphi}e^{i\mu\varphi}\right)\left(\frac{d}{d\varphi}\cos(\varphi)^\alpha\right) d\varphi=\\
&\left(\frac{d}{d\varphi}e^{i\mu\varphi}\right)\cos(\varphi)^\alpha |_{-\beta}^{\beta}-e^{i\mu\varphi}\left(\frac{d}{d\varphi}\cos(\varphi)^\alpha\right) |_{-\beta}^{\beta}+\int_{-\beta}^{\beta}e^{i\mu\varphi}\left(\frac{d^2}{d\varphi^2}\cos(\varphi)^\alpha\right) d\varphi=\\
&\Big(e^{i\mu\varphi}\cos(\varphi)^{\alpha-1} (i\mu \cos(\varphi)+\alpha \sin(\varphi))\Big) \Big|_{-\beta}^{\beta}+\int_{-\beta}^{\beta}e^{i\mu\varphi}\left(\frac{d^2}{d\varphi^2}\cos(\varphi)^\alpha\right) d\varphi=\\
&\Big(e^{i\mu\varphi}\cos(\varphi)^{\alpha-1} (i\mu \cos(\varphi)+\alpha \sin(\varphi))\Big) \Big|_{-\beta}^{\beta}+\int_{-\beta}^{\beta}e^{i\mu\varphi}\Big(\alpha(\alpha-1)\cos(\varphi)^{\alpha-2}-\alpha^2\cos(\varphi)^\alpha\Big)d\varphi=\\
&
\Big(e^{i\mu\varphi}\cos(\varphi)^{\alpha-1} (i\mu \cos(\varphi)+\alpha \sin(\varphi))\Big) \Big|_{-\beta}^{\beta}+\alpha(\alpha-1)P(\mu,\alpha-2,\beta)-\alpha^2P(\mu,\alpha,\beta).
\end{align*}
Hence, we obtain the following recursive relation for $P$
\begin{align*}
&P(\mu,\alpha-2,\beta)=\frac{\alpha}{\alpha-1}\left(1-\frac{\mu^2}{\alpha^2}\right)P(\mu,\alpha,\beta)+h(\mu,\alpha,\beta),
\end{align*}
where
\begin{align}\label{hdefin}
&h(\mu,\alpha,\beta)=
\frac{2}{\alpha(\alpha-1)}(\mu\cos(\beta)^\alpha\sin(\mu\beta)-\alpha \cos^{\alpha-1}(\beta)\cos(\mu\beta)\sin(\beta)).
\end{align}
For all positive $n$ it thus holds that
\begin{align}
\nonumber
&P(\mu,\alpha-2,\beta)=\prod_{k=0}^{n}\left( \frac{\alpha+2k}{\alpha+2k-1} \right)
\prod_{k=0}^{n}\left( (1-\frac{\mu^2}{(\alpha+2k)^2} \right)P(\mu,\alpha+2n,\beta)+\\ \nonumber&\sum_{j=0}^{n}\prod_{k=0}^{j-1}\left( \frac{\alpha+2k}{\alpha+2k-1} \right)
\prod_{k=0}^{j-1}\left( (1-\frac{\mu^2}{(\alpha+2k)^2} \right)h(\mu,\alpha+2j,\beta)=\\
\nonumber&\frac{\Gamma(\frac{\alpha-1}{2})}{\Gamma(\frac{\alpha-1}{2}+n)} \frac{\Gamma(\frac{\alpha}{2}+n)}{\Gamma(\frac{\alpha}{2})}
\prod_{k=0}^{n}\left( (1-\frac{\mu^2}{(\alpha+2k)^2} \right)P(\mu,\alpha+2n,\beta)+\\ \nonumber&\sum_{j=0}^{n}\frac{\Gamma(\frac{\alpha-1}{2})}{\Gamma(\frac{\alpha-1}{2}+j-1)} \frac{\Gamma(\frac{\alpha}{2}+j-1)}{\Gamma(\frac{\alpha}{2})}
\prod_{k=0}^{j-1}\left( (1-\frac{\mu^2}{(\alpha+2k)^2} \right)h(\mu,\alpha+2j,\beta)=\\
\nonumber&\frac{\mathcal{B}\left(\frac{\alpha}{2}-\frac{1}{2},\frac{1}{2}\right)}{\mathcal{B}(n+\frac{\alpha}{2}-\frac{1}{2},\frac{1}{2})}
\prod_{k=0}^{n}\left( (1-\frac{\mu^2}{(\alpha+2k)^2} \right)P(\mu,\alpha+2n,\beta)+\\ &\label{expr1P}\sum_{j=0}^{n}\frac{\mathcal{B}(\frac{\alpha}{2}-\frac{1}{2},\frac{1}{2})}{\mathcal{B}(j+\frac{\alpha}{2}-\frac{1}{2},\frac{1}{2})}
\prod_{k=0}^{j-1}\left( (1-\frac{\mu^2}{(\alpha+2k)^2} \right)h(\mu,\alpha+2j,\beta),
\end{align}
where $\mathcal{B}$ is the beta function, defined by
\begin{equation}\label{betadef}
\mathcal{B}(m,n)=2\int_0^{\pi/2}(\cos \varphi)^{2m-1}(\sin \varphi)^{2n-1}d\varphi.
\end{equation}
Above we have used some of the properties of the gamma function \cite[p. 256]{abramowitz1972handbook} and its relation to the beta function $\mathcal{B}$  \cite[p. 258]{abramowitz1972handbook},
\begin{align}\label{betaprop}
\mathcal{B}(m,n)=\frac{\Gamma(m)\Gamma(n)}{\Gamma(m+n)}.
\end{align}
In the limit case $n\to\infty$, it holds that
\begin{align}
\nonumber&P(\mu,\alpha-2,\beta)=\lim_{n\to\infty}\Bigg[
\frac{\mathcal{B}(\frac{\alpha}{2}-\frac{1}{2},\frac{1}{2})}{\mathcal{B}(n+\frac{\alpha}{2}-\frac{1}{2},\frac{1}{2})}
\prod_{k=0}^{n}\left( 1-\frac{\mu^2}{(\alpha+2k)^2} \right)P(\mu,\alpha+2n,\beta)+\\ \nonumber&\sum_{j=0}^{n}\frac{\mathcal{B}(\frac{\alpha}{2}-\frac{1}{2},\frac{1}{2})}{\mathcal{B}(j+\frac{\alpha}{2}-\frac{1}{2},\frac{1}{2})}
\prod_{k=0}^{j-1}\left( 1-\frac{\mu^2}{(\alpha+2k)^2} \right)h(\mu,\alpha+2j,\beta)\Bigg]=\\
\nonumber&\lim_{n\to\infty}\Bigg[
\mathcal{B}\left(\frac{\alpha}{2}-\frac{1}{2},\frac{1}{2}\right)
\prod_{k=0}^{n}\left( 1-\frac{\mu^2}{(\alpha+2k)^2} \right)\int_{-\beta}^{\beta}e^{i\mu\varphi}\frac{\cos(\varphi)^{\alpha+2n}}{\mathcal{B}(n+\frac{\alpha}{2}-\frac{1}{2},\frac{1}{2})} d\varphi+\\ \label{exprlim}&\sum_{j=0}^{n}\frac{\mathcal{B}(\frac{\alpha}{2}-\frac{1}{2},\frac{1}{2})}{\mathcal{B}(j+\frac{\alpha}{2}-\frac{1}{2},\frac{1}{2})}
\prod_{k=0}^{j-1}\left( 1-\frac{\mu^2}{(\alpha+2k)^2} \right)h(\mu,\alpha+2j,\beta)\Bigg]
\end{align}
Using \eqref{betadef}, it is easy to see that the sequence of functions
\begin{align*}
h_n(t)=\frac{(\cos t)^n}{\mathcal{B}(\frac{n}{2}-\frac{1}{2},\frac{1}{2})}
\end{align*}
tends to $\delta(t)$ when $n\to\infty$. Using this fact and \eqref{exprlim}
we can represent $P(\mu,\alpha,\beta)$ in the form $P(\mu,\alpha,\beta)=P_0(\mu,\alpha,\beta)+P_1(\mu,\alpha,\beta)$, where
\begin{align*}
&P_0(\mu,\alpha,\beta)=\mathcal{B}\left(\frac{\alpha}{2}+\frac{1}{2},\frac{1}{2}\right)
\prod_{k=1}^{\infty}\left( 1-\frac{\mu^2}{(\alpha+2k)^2} \right),\\
&P_1(\mu,\alpha,\beta)=\sum_{j=1}^{\infty}\frac{\mathcal{B}(\frac{\alpha}{2}+\frac{1}{2},\frac{1}{2})}{\mathcal{B}(j+\frac{\alpha}{2}+\frac{1}{2},\frac{1}{2})}
\prod_{k=1}^{j-1}\left( 1-\frac{\mu^2}{(\alpha+2k)^2} \right)h(\mu,\alpha+2j,\beta).
\end{align*}
From the identity \cite[p. 336]{ramanujan1994P4}
\begin{align}
\label{identity}
\frac{\Gamma^2(n+1)}{\Gamma(n+x i+1)\Gamma(n-x i+1)}=\prod_{k=1}^{\infty}\left(1+\frac{x^2}{(n+k)^2}\right),
\end{align}
it follows that
\begin{align}
P_0(\mu,\alpha,\beta) &= \mathcal{B}\left(\frac{\alpha+1}{2},\frac{1}{2}\right)
\prod_{k=1}^{\infty}\left( 1-\frac{\mu^2}{(\alpha+2k)^2} \right) \label{exprgamma}\\
&=\mathcal{B}\left(\frac{\alpha+1}{2},\frac{1}{2}\right)\frac{\Gamma^2(\frac{\alpha}{2}+1)}{\Gamma(\frac{\alpha}{2}+\frac{\mu}{2}+1)\Gamma(\frac{\alpha}{2}-\frac{\mu}{2}+1)}=\frac{\Gamma(\frac{\alpha+1}{2})\Gamma(\frac{1}{2})\Gamma(\frac{\alpha+2}{2})}{\Gamma(\frac{\alpha+\mu}{2}+1)\Gamma(\frac{\alpha-\mu}{2}+1)}. \nonumber
\end{align}
Concerning  the second term, we have that
\begin{align*}
&P_1(\mu,\alpha,\beta)= \lim_{n\to\infty}\left[
\sum_{j=0}^{n}\frac{\mathcal{B}(\frac{\alpha+1}{2},\frac{1}{2})}{\mathcal{B}(j+\frac{\alpha+1}{2},\frac{1}{2})}h(\mu,\alpha+2+2j)
\prod_{k=1}^{j}\left( 1-\frac{\mu^2}{(\alpha+2k)^2} \right)\right],
\end{align*}
where $h(\mu,\alpha )=
\frac{2}{\alpha(\alpha-1)}(\mu\cos(\beta)^\alpha\sin(\mu\beta)-\alpha \cos^{\alpha-1}(\beta)\cos(\mu\beta)\sin(\beta))$. By using the identity \eqref{identity} again, we can rewrite
\begin{align*}
&\prod_{k=1}^{j}\left( 1-\frac{\mu^2}{(\alpha+2k)^2} \right)=
\frac{\Gamma^2(\frac{\alpha}{2}+1)\Gamma(\frac{\alpha}{2}+\frac{\mu}{2}+j+1)\Gamma(\frac{\alpha}{2}-\frac{\mu}{2}+j+1)}{\Gamma(\frac{\alpha}{2}+\frac{\mu}{2}+1)\Gamma(\frac{\alpha}{2}-\frac{\mu}{2}+1)\Gamma^2(\frac{\alpha}{2}+j+1)}
\end{align*}
An expression for $P_1(\mu,\alpha,\beta)$ thus reads,
\begin{align}
\nonumber P_1(\mu,\alpha,\beta)&=
\sum_{j=0}^{\infty}
h(\mu,\alpha+2+2j)
\frac{\Gamma\left(\frac{\alpha+1}{2}\right)\Gamma(\frac{\alpha}{2}+1)\Gamma(\frac{\alpha}{2}+\frac{\mu}{2}+j+1)\Gamma(\frac{\alpha}{2}-\frac{\mu}{2}+j+1)}{\Gamma(j+\frac{\alpha+1}{2})\Gamma(\frac{\alpha}{2}+\frac{\mu}{2}+1)\Gamma(\frac{\alpha}{2}-\frac{\mu}{2}+1)\Gamma(\frac{\alpha}{2}+j+1)}\\
&=\sum_{j=0}^{\infty}
h(\mu,\alpha+2+2j)
\frac{(\frac{\alpha}{2}+\frac{\mu}{2}+1)_j(\frac{\alpha}{2}-\frac{\mu}{2}+1)_j}{(\frac{\alpha+1}{2})_j(\frac{\alpha}{2}+1)_j}\\
\end{align}
using the  Pochhammer symbol \cite[p. 256]{abramowitz1972handbook} 
\begin{align*}
&(x)_j=\frac{\Gamma(x+j)}{\Gamma(x)}=x(x+1)\dots(x+j-1).
\end{align*}
By using the definition of $h$ \eqref{hdefin}, we split the sum into two parts,
\begin{align*}
&\nonumber P_1(\mu,\alpha,\beta)=\sum_{j=0}^{\infty}
\frac{2\mu\cos(\beta)^{\alpha+2+2j}\sin(\mu\beta)}{(\alpha+2+2j)(\alpha+1+2j)}
\frac{(\frac{\alpha}{2}+\frac{\mu}{2}+1)_j(\frac{\alpha}{2}-\frac{\mu}{2}+1)_j}{(\frac{\alpha+1}{2})_j(\frac{\alpha}{2}+1)_j}-\\
&\sum_{j=0}^{\infty}
\frac{2 \cos^{\alpha+1+2j}(\beta)\cos(\mu\beta)\sin(\beta)}{(\alpha+1+2j)}
\frac{(\frac{\alpha}{2}+\frac{\mu}{2}+1)_j(\frac{\alpha}{2}-\frac{\mu}{2}+1)_j}{(\frac{\alpha+1}{2})_j(\frac{\alpha}{2}+1)_j}=\\
&2\mu\cos(\beta)^{\alpha+2}\sin(\mu\beta)\sum_{j=0}^{\infty}
\frac{\cos(\beta)^{2j}}{(\alpha+2+2j)(\alpha+1+2j)}
\frac{(\frac{\alpha}{2}+\frac{\mu}{2}+1)_j(\frac{\alpha}{2}-\frac{\mu}{2}+1)_j}{(\frac{\alpha+1}{2})_j(\frac{\alpha}{2}+1)_j}-\\
&2\cos(\beta)^{\alpha+1}\cos(\beta \mu)\sin(\beta)\sum_{j=0}^{\infty}
\frac{\cos(\beta)^{2j}}{(\alpha+1+2j)}
\frac{(\frac{\alpha}{2}+\frac{\mu}{2}+1)_j(\frac{\alpha}{2}-\frac{\mu}{2}+1)_j}{(\frac{\alpha+1}{2})_j(\frac{\alpha}{2}+1)_j}.
\end{align*}
Now we can use the facts that $(1)_j=j!$ and that
\begin{align*}
(z+j)(z)_j=(z)_{j+1}=\frac{\Gamma(z+j+1)}{\Gamma(z+1)}=\frac{1}{z}\frac{\Gamma(z+j+1)}{\Gamma(z)}=\frac{1}{z}(z+1)_j,
\end{align*}
to simplify the above expression as
\begin{align*}
&P_1(\mu,\alpha,\beta)=
2\mu\cos(\beta)^{\alpha+2}\sin(\mu\beta)\frac{1}{(\alpha+1)(\alpha+2)}\sum_{j=0}^{\infty}
\frac{\cos(\beta)^{2j}}{j!}
\frac{(1)_j(\frac{\alpha}{2}+\frac{\mu}{2}+1)_j(\frac{\alpha}{2}-\frac{\mu}{2}+1)_j}{(\frac{\alpha+3}{2})_j(\frac{\alpha}{2}+2)_j}-\\
&2\cos(\beta)^{\alpha+1}\cos(\beta \mu)\sin(\beta)\frac{1}{(\alpha+1)}\sum_{j=0}^{\infty}
\frac{\cos(\beta)^{2j}}{j!}
\frac{(1)_j(\frac{\alpha}{2}+\frac{\mu}{2}+1)_j(\frac{\alpha}{2}-\frac{\mu}{2}+1)_j}{(\frac{\alpha+3}{2})_j(\frac{\alpha}{2}+1)_j}.
\end{align*}
The sums above take the form of $_3F_2$ hypergeometric functions \cite[p. 403]{olver2010nist}, and hence it holds that
\begin{align*}
P_1(\mu,\alpha,\beta) &= \frac{2\mu\cos(\beta)^{\alpha+2}\sin(\mu\beta)}{(\alpha+1)(\alpha+2)}\,_3F_2\left(1,\frac{\alpha}{2}+\frac{\mu}{2}+1,\frac{\alpha}{2}-\frac{\mu}{2}+1;\frac{\alpha+3}{2},\frac{\alpha}{2}+2;\cos(\beta)^2\right)\\
&-\frac{2\cos(\beta)^{\alpha+1}\cos(\beta \mu)\sin(\beta)}{(\alpha+1)} \,_3F_2\left(1,\frac{\alpha}{2}+\frac{\mu}{2}+1,\frac{\alpha}{2}-\frac{\mu}{2}+1;\frac{\alpha+3}{2},\frac{\alpha}{2}+1;\cos(\beta)^2\right).
\end{align*}
We then finally obtain that
\begin{align*}
&P(\mu,\alpha,\beta)=
\frac{\Gamma(\frac{\alpha+1}{2})\Gamma(\frac{1}{2})\Gamma(\frac{\alpha+2}{2})}{\Gamma(\frac{\alpha+\mu}{2}+1)\Gamma(\frac{\alpha-\mu}{2}+1)}+\\
&\frac{2\mu\cos(\beta)^{\alpha+2}\sin(\mu\beta)}{(\alpha+1)(\alpha+2)}\,_3F_2\left(1,\frac{\alpha}{2}+\frac{\mu}{2}+1,\frac{\alpha}{2}-\frac{\mu}{2}+1;\frac{\alpha+3}{2},\frac{\alpha}{2}+2;\cos(\beta)^2\right)-\\
&\frac{2\cos(\beta)^{\alpha+1}\cos(\beta \mu)\sin(\beta)}{(\alpha+1)} \,_3F_2\left(1,\frac{\alpha}{2}+\frac{\mu}{2}+1,\frac{\alpha}{2}-\frac{\mu}{2}+1;\frac{\alpha+3}{2},\frac{\alpha}{2}+1;\cos(\beta)^2\right).
\end{align*}
\bibliographystyle{siam}
\bibliography{lp_ref}
\end{document}